\documentclass[final,sort&compress,12pt]{elsarticle}
\pdfoutput=1
\usepackage{fullpage} 
\usepackage[inline,shortlabels]{enumitem}   
\setlist[enumerate,itemize]{topsep=2pt,itemsep=-0.5ex,partopsep=1ex,parsep=1ex}
\usepackage{amsmath}	        	
\usepackage{amsthm}							
\usepackage{amssymb}			
\usepackage{mathtools}					
\usepackage{mathrsfs}						
\usepackage{extarrows}          
\usepackage{xfrac}              
\usepackage{url} 
\usepackage[normalem]{ulem}

\usepackage{anyfontsize}
\PassOptionsToPackage{hyphens}{url}
\usepackage[hidelinks,unicode,verbose]{hyperref}
\usepackage{color}
\usepackage{caption}
\usepackage{subcaption}
\usepackage{tocbibind}
\usepackage{txfonts}
\journal{Stochastic Processes and Their Applications}

\makeatletter
\def\ps@pprintTitle{%
     \let\@oddhead\@empty
		 \let\@evenhead\@empty
     \def\@oddfoot
       {\hbox to \textwidth%
        {\ifnopreprintline\relax\else
        \@myfooterfont%
         \ifx\@elsarticlemyfooteralign\@elsarticlemyfooteraligncenter%
           \hfil\@elsarticlemyfooter\hfil%
         \else%
         \ifx\@elsarticlemyfooteralign\@elsarticlemyfooteralignleft%
           \@elsarticlemyfooter\hfill{}%
         \else%
         \ifx\@elsarticlemyfooteralign\@elsarticlemyfooteralignright%
           {}\hfill\@elsarticlemyfooter%
         \else%
            \begin{minipage}[t]{\textwidth}   Preprint accepted in \ifx\@journal\@empty%
                 Elsevier%
            \else\@journal\fi\hfill March 28, 2023\\
						Published version available at \rurl{10.1016/j.spa.2023.04.010}\medskip\\ 
						\begin{minipage}[b]{0.7\textwidth}\textcopyright\ 2023. This work is licensed under a Creative Commons\\
			      Attribution-NonCommercial-NoDerivatives
						4.0 International License\end{minipage} 
							\hfill\includegraphics[width=0.2\textwidth]{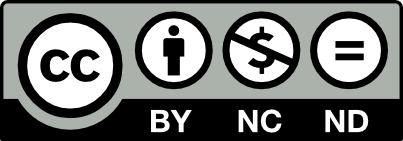}
						\end{minipage}\fi%
         \fi%
         \fi%
         \fi%
         }
       }%
     \let\@evenfoot\@oddfoot}
\makeatother
\newcommand\rurl[1]{%
  \href{http://doi.org/#1}{\nolinkurl{#1}}%
}

\let\originalleft\left
\let\originalright\right
\renewcommand{\left}{\mathopen{}\mathclose\bgroup\originalleft}
\renewcommand{\right}{\aftergroup\egroup\originalright}

\makeatletter
\newcommand*{\bigcdot}{}
\DeclareRobustCommand*{\bigcdot}{%
  \mathbin{\mathpalette\bigcdot@{}}%
}
\newcommand*{\bigcdot@scalefactor}{.7}
\newcommand*{\bigcdot@widthfactor}{1.15}
\newcommand*{\bigcdot@}[2]{%
  \sbox0{$#1\vcenter{}$}
  \sbox2{$#1\cdot\m@th$}%
  \hbox to \bigcdot@widthfactor\wd2{%
    \hfil
    \raise\ht0\hbox{%
      \scalebox{\bigcdot@scalefactor}{%
        \lower\ht0\hbox{$#1\bullet\m@th$}%
      }%
    }%
    \hfil
  }%
}
\makeatother

\newcommand{\widebar}[1]{\overbracket[0.5pt][0pt]{#1}}
\newcommand{\underli}[1]{\mkern2mu\mbox{\uline{$\mkern-2mu#1\mkern-4mu$}}\mkern4mu }

\newcommand{\N}{\mathbb{N}}

\newcommand{\R}{\mathbb{R}}
\newcommand{\Cx}{\mathbb{C}}
\newcommand{\Cinf}{{\widebar{\Cx}\mkern0.5mu}}
\newcommand{\Oinf}{{\widebar{\Omega}\mkern0.5mu}}

\newcommand{\Uni}{\mathfrak{U}}
\newcommand{\I}{\mathfrak{I}}
\renewcommand{\U}{\mathcal{U}}

\renewcommand{\d}{\mathrm{d}}
\newcommand{\id}{\mathord{\rm id}}
\newcommand{\e}{\mathrm{e}}
\renewcommand{\Re}{\operatorname{Re}}
\renewcommand{\Im}{\operatorname{Im}}
\DeclareMathOperator{\sgn}{sgn}

\DeclareMathOperator{\E}{{}\textsf{\upshape E}}

\renewcommand{\P}{\textsf{\upshape P}}
\newcommand{\Qu}{\textsf{\upshape Q}}
\newcommand{\indicator}[1]{\mathbf{1}_{#1}}

\newcommand{\filt}[1]{\mathfrak{#1}}
\newcommand{\sigalg}[1]{\mathscr{#1}}

\newcommand{\T}{\mathcal{T}}
\newcommand{\ddp}{\mathrm{dp}}
\newcommand{\qc}{\mathrm{qc}}
\newcommand{\Exp}{\ensuremath{\mathscr{E}}}
\newcommand{\Log}{\ensuremath{\mathcal{L}}}
\newcommand{\V}{\ensuremath{\mathscr{V}}}
\newcommand{\sint}{\bigcdot}

\newcommand{\NaN}{\mathrm{NaN}}
\newcommand{\lc}{[\![}
\newcommand{\rc}{]\!]}

\newcommand*{\bigs}[1]{\scalebox{1.2}{\ensuremath#1}}
\newcommand\bigsl[1]{\mathopen{\bigs{#1}}}
\newcommand\bigsr[1]{\mathclose{\bigs{#1}}}
\DeclarePairedDelimiter\abs{|}{|}

\theoremstyle{plain}
\newtheorem{theorem}{Theorem}
\newtheorem{lemma}[theorem]{Lemma}
\newtheorem{proposition}[theorem]{Proposition}
\newtheorem{corollary}[theorem]{Corollary}
\theoremstyle{definition}
\newtheorem{definition}[theorem]{Definition}

\newtheorem{example}[theorem]{Example}

\theoremstyle{remark}
\newtheorem{remark}[theorem]{Remark}

\numberwithin{theorem}{section}
\numberwithin{equation}{section}
\numberwithin{figure}{section}
\numberwithin{table}{section}

\let\OLDthebibliography\thebibliography
\renewcommand\thebibliography[1]{
  \OLDthebibliography{#1}
  \setlength{\parskip}{0pt}
  \setlength{\itemsep}{1pt plus 0.3ex}
}

\begin{document}

\begin{frontmatter}

\title{Simplified Calculus for Semimartingales:\\ {M}ultiplicative Compensators and Changes of Measure}

\author[bayes]{Ale\v{s} \v{C}ern\'{y}\texorpdfstring{\,\corref{cor1}}{}}
\ead{ales.cerny.1@city.ac.uk}
\author[lse]{Johannes Ruf\texorpdfstring{\,}{}}
\ead{j.ruf@lse.ac.uk}

\cortext[cor1]{Corresponding author}
\address[bayes]{Bayes Business School, City, University of London, 106 Bunhill Row, London EC1Y 8TZ, UK}
\address[lse]{Department of Mathematics, London School of Economics and Political Science, Columbia House, Houghton Street, London WC2A 2AE, UK}
\begin{abstract}
The paper develops multiplicative compensation for complex-valued semimartingales and studies some of its consequences. It is shown that the stochastic exponential of any complex-valued semimartingale with independent increments becomes a true martingale after multiplicative compensation when such compensation is meaningful. This generalization of the L\'evy--Khintchin formula fills an existing gap in the literature. It allows, for example, the computation of the Mellin transform of a signed stochastic exponential, which in turn has practical applications in mean--variance portfolio theory. Girsanov-type results based on multiplicatively compensated semimartingales simplify treatment of absolutely continuous measure changes. As an example, we obtain the characteristic function of log returns for a popular class of minimax measures in a L\'evy setting.
\end{abstract}

\end{frontmatter}
\setlength{\abovedisplayskip}{4pt}
\setlength{\belowdisplayskip}{4pt}

\section{Introduction}\label{sect: intro}
Multiplicatively compensated semimartingales are an important tool in stochastic modelling. They appear, among others, in the following contexts.
\begin{itemize}
\item \emph{Computation of characteristic functions}, e.g., Jacod and Shiryaev \cite[Section~III.7]{js.03}.
\item \emph{Esscher-type measure changes} including a variety of minimax martingale measures, 
e.g., Goll and R\"uschendorf \cite{goll.rueschendorf.01}, Jeanblanc et al. \cite{jeanblanc.al.07}; martingale measures associated with ad-hoc numeraire changes, e.g., Eberlein et al. \cite{eberlein.al.09}; 
but also non-martingale measures, e.g., the opportunity-neutral measure in semimartingale mean-variance theory, e.g., \v{C}ern\'y and Kallsen \cite[Section~3.4]{cerny.kallsen.07}. A very general take on Esscher-type measures is presented in Kallsen and Shiryaev \cite{kallsen.shiryaev.02}.
\item \emph{Proofs of moment bounds}, e.g., to show existence and uniqueness of BSDE solutions, e.g., Kazi-Tani et al. \cite[Lemma~A.5]{kazi-tani.al.15}; to estimate variation distance of probability measures in Kabanov et al. \cite[Theorem~2.1]{kabanov.al.86}; to prove uniform integrability of local martingales, e.g., L\'epingle and M\'emin \cite[Th\'eor\`eme~III.1]{lepingle.memin.78.ptrf}, Ruf \cite[Corollary~5]{Ruf_Novikov}.
\item \emph{Filtration extension and/or shrinkage}, e.g., Nikeghbali and Yor \cite[Section~4]{nikeghbali.yor.06}, Kardaras \cite{Kardaras_times},  Aksamit and Jeanblanc \cite[Chapter~5]{aksamit.jeanblanc.17};  Kardaras and Ruf \cite[Section~5]{kardaras.ruf.20}.
\item \emph{Theory of Markov processes}, e.g., It\^o and Watanabe \cite[Chapter~2]{ito.watanabe.65}; Chen et al. \cite[Theorem~3.1]{chen.al.04}.
\end{itemize}

This paper examines some of the consequences of multiplicative compensation for signed (and even complex-valued) semimartingales. The next statement is a special case of Theorem~\ref{T:200603}\ref{T:200603.1}.
\begin{theorem}\label{T:intro1}
Let $Y$ be a special $\Cx$--valued semimartingale with independent increments. Then 
$$\E[\Exp(Y)_t]=\Exp\bigs(B^Y\bigs)_t,\qquad t>0.$$
\end{theorem}  
\noindent Here $\Exp(Y)$ denotes the stochastic exponential%
\footnote{By convention, any stochastic exponential starts at $1$, i.e., here $\Exp(Y)_0=1$.}
of a semimartingale $Y$ and $B^Y$ the predictable finite-variation part (here also called the drift) in the canonical decomposition of a special semimartingale~$Y$.%
\footnote{We assume that $B^Y_0=0$, which makes the finite variation part $B^Y$ unique; see \cite[I.4.22]{js.03}.}
Theorem~\ref{T:intro1} with $Y = (\e^{iu\hspace{0.05em}\id}-1)\circ\mkern-2mu X$, $X$ a L\'evy process, $u\in\R$, and $\id$ the identity function  recovers the L\'evy--Khintchin formula for $X$ (Corollary~\ref{C:LKformula}).  The operation $f\circ\mkern-2mu X$ denotes, roughly speaking, the $f$--variation of $X$ (Definition~\ref{D:181101}). Observe that in Theorem~\ref{T:intro1} the process $\Exp(B^Y)$  is deterministic.

To illustrate the novelty of Theorem~\ref{T:intro1}, consider the task of computing the distribution of $\Exp(X)$ when the stochastic exponential is signed. Here, one can evaluate $|\Exp(X)|^{iu}$ and $\sgn(\Exp(X))|\Exp(X)|^{iu}$ for $u\in\R$ separately to obtain the Mellin transforms of the positive and negative parts, respectively (Example~\ref{E:190218}),
\begin{align}
\E\bigsl[|\Exp(X)_t|^{iu}\bigsr] ={}& \Exp\bigs(B^{(|1+\id|^{iu} - 1)\circ\mkern-1mu X}\bigs)_t;\label{eq:210404.4}\\
\E\bigsl[\sgn(\Exp(X)_t)|\Exp(X)_t|^{iu}\bigsr] ={}& \Exp\bigs(B^{(\sgn(1+\id)|1+\id|^{iu} - 1)\circ\mkern-1mu X}\bigs)_t.\label{eq:210404.5}
\end{align}
The standard L\'evy-Khintchin formula is unable to deliver such a result. The Mellin transform can then be used  to solve, at least numerically, previously intractable  questions concerning mean-variance portfolio allocation (Example~\ref{E:MV}).

Next, as a special case of Theorem~\ref{T:mult comp}\ref{MC:2} below, we have a complex-valued extension of a classical but perhaps not sufficiently well-known result due to L\'epingle and M\'emin \cite[Proposition~II.1]{lepingle.memin.78.ptrf}.
\begin{theorem}\label{T:intro2}
Let $Y$ be a special $\Cx$--valued semimartingale such that $\Delta B^Y \neq -1$. Then $\Exp(B^Y)$ is a multiplicative compensator of $\Exp(Y)$, that is, $\frac{\Exp(Y)}{\Exp(B^Y)}$ is a local martingale.
\end{theorem}
\noindent  In the main body of the paper (see Theorem~\ref{T:mult comp}), we allow $Y$ to explode on approach to a stopping time, which later allows study of a larger class of non-equivalent measure changes. In the context of Theorem~\ref{T:intro2}, Theorem~\ref{T:intro1} asserts that $\frac{\Exp(Y)}{\Exp(B^Y)}$ is a true martingale whenever $Y$ has independent increments (Theorem~\ref{T:200603}\ref{T:200603.2}).

Multiplicative compensation of semimartingales is, of course, not new (see, e.g., Yoeurp and Meyer \cite{yoeurp.meyer.76}, Az\'ema \cite{azema.78}, Jacod \cite{jacod.78}, M\'emin \cite{memin.78}, L\'epingle and M\'emin \cite{lepingle.memin.78.ptrf}, and Kallsen and Shiryaev \cite{kallsen.shiryaev.02}).  We wish to emphasize the strength and flexibility of Theorem~\ref{T:intro2} when coupled with semimartingale representations. For example, Theorem~2.19 in \cite{kallsen.shiryaev.02} follows immediately from Theorem~\ref{T:intro2} by writing (see Example~\ref{E:190703.1})
$$\Exp(Y) = \e^{{\zeta}\sint{X}} = \Exp\bigsl(\bigsl(\e^{\zeta\id}-1\bigsr)\circ\mkern-2mu X\bigsr),\qquad \zeta\in L(X),$$
where $\zeta\sint X$ stands for the stochastic integral $\int_0^\cdot \zeta_t \d X_t$. Many computations in \cite{lepingle.memin.78.ptrf} follow from Theorem~\ref{T:intro2} by taking (for $\lambda > 0$ and  a local martingale $M$),
$$\Exp(Y) = \Exp^\lambda(M) = \Exp\bigsl(\bigsl(\bigsl(1+\id\bigsr)^\lambda - 1\bigsr)\circ M\bigsr),\qquad \qquad  \text{etc.}$$

Theorem~\ref{T:intro2} is a stepping stone to Girsanov-type results for measure changes relying on non-negative multiplicatively compensated semimartingales.
\begin{theorem}\label{T:intro3}
Let $Y$ be a special $\R$--valued semimartingale such that $\Delta Y > -1$ 
and ${M=\frac{\Exp(Y)}{\Exp(B^Y)}}$ is a uniformly integrable martingale. For the measure $\Qu$ given by 
$\frac{\d\Qu}{\d\P}=M_\infty$ and a semimartingale $X$, the following are equivalent.
\begin{enumerate}[label={\rm(\roman{*})}, ref={\rm(\roman{*})}]
\item\label{T:intro3:i} $X$ is $\Qu$--special.
\item\label{T:intro3:ii} $X\frac{\Exp(Y)}{\Exp(B^Y)}$ is $\P$--special.
\item\label{T:intro3:iii} $X +[X,Y]$ is $\P$--special.
\end{enumerate}
If either of these conditions holds, then $B^{X}_\Qu=\frac{1}{1+\Delta B^Y}\sint B^{X +[X,Y]}$ and
$$\Exp(B^{X}_\Qu) = \frac{\Exp\bigs(B^{X+Y+[X,Y]}\bigs)}{\Exp(B^Y)},$$
where $B^{X}_\Qu$ stands for the additive $\Qu$--compensator of the $\Qu$--special semimartingale $X$. 
Furthermore, the following are equivalent.
\begin{enumerate}[label={\rm(\roman{*}')}, ref={\rm(\roman{*}')}]
\item\label{T:intro3:i'} $X$ is a $\Qu$--local martingale.
\item\label{T:intro3:ii'} $X\frac{\Exp(Y)}{\Exp(B^Y)}$ is a $\P$--local martingale.
\item\label{T:intro3:iii'} $X +[X,Y]$ is a $\P$--local martingale.
\end{enumerate} 
\end{theorem}
\begin{proof}
See Theorem~\ref{T:Girsanov} and Corollary~\ref{C:181005} below.
\end{proof}

As a practical application of Theorem~\ref{T:intro3}, Example~\ref{E:210410} in the main body of the paper computes the characteristic function of log returns for a popular class of minimax measures. Additionally, Example~\ref{E:221123} showcases the usefulness of Theorem~\ref{T:intro3} by streamlining an otherwise fairly involved calculation appearing in \cite{lepingle.memin.78.ptrf}. 

In the rest of the paper we proceed as follows. Section~\ref{S:setup} provides the setup of this paper.  Section~\ref{S:3} discusses the construction of multiplicative compensators and provides several examples. Section~\ref{S:4} considers multiplicative compensation for stochastic exponentials of processes with independent increments. Section~\ref{S:5} introduces a version of Girsanov's theorem and considers additive and multiplicative compensation after a change of measure. Finally, Section~\ref{S:6} concludes.

\section{Setup and notation}\label{S:setup}
The applications in this paper rely on semimartingale representations worked out in \v{C}ern\'y and Ruf \cite{cerny.ruf.22.ejp}. A whittled-down summary of the relevant definitions and results from \cite{cerny.ruf.22.ejp} is provided in Subsections~\ref{SS:2.1}--\ref{SS:2.5}. We suggest skipping the details on the first reading and instead making use of the following ``executive summary.'' One can also consult the introductory paper \cite{cerny.ruf.21.ejor} for further context and examples.

For a  function $\xi:[0,\infty)\times \R\to\R$ that is constant in the first argument and twice continuously differentiable in the second argument and for an $\R$--valued semimartingale $X$, the partial sums 
$$ \sum_{n\in\N}\xi_{t_{n-1}} \left(X^{t_n}-X^{t_{n-1}}\right)$$ 
converge uniformly on compact time intervals in probability to 
\begin{equation}\label{eq:emery}
\xi\circ\mkern-2mu X := \xi'(0)\sint X+\frac{1}{2}\xi''(0)\sint [X,X]^c + \sum_{0<t\leq\cdot}(\xi_t(\Delta X_t)-\xi_t'(0)\Delta X_t)
\end{equation}
as the time partition $(t_n)_{n\in\N}$ becomes finer; see \'Emery~\cite[Th\'eor\`eme~2a]{emery.78}. In such case, $\xi \circ\mkern-2mu X$ can be interpreted as the $\xi$--variation of $X$. Here and below we write $X^\tau$ for the process $X$ stopped at some stopping time $\tau$.

Formula \eqref{eq:emery} also makes sense  for some predictable functions $\xi$. In particular, it makes sense for all predictable functions in the universal class $\Uni$ defined below, including all twice-continuously differentiable deterministic time-constant functions. For $\xi\in\Uni$, $\xi\circ\mkern-2mu X$ is not necessarily a $\xi$--variation. We then speak more broadly of semimartingale representations, saying that $Y$ is represented in terms of $X$ if there is $\xi\in\Uni$ such that $Y=\xi\circ\mkern-2mu X$. The functions in $\Uni$, such as $\log(1+\id)$, need not be defined everywhere. For $\xi\circ\mkern-2mu X$ to exist, it is enough that $X$ is compatible with $\xi\in\Uni$, i.e., $\xi(\Delta X)$ is finite almost surely.

The \'Emery formula \eqref{eq:emery} has a natural counterpart for complex-valued functions of several complex variables and the definition of $\Uni$ takes this into account. Semimartingale representations in $\Uni$ conveniently capture common operations on semimartingales. For example, locally bounded stochastic integration corresponds to ``linear variation'' 
$$\zeta\sint X = (\zeta\id) \circ\mkern-2mu X;$$ 
smooth transformation, too, has a simple representation in $\Uni$,
$$f(X)-f(X_0) = (f(X_-+\id)-f(X_-))\circ\mkern-2mu X.$$ 
Furthermore, $\Uni$ is closed under composition, with 
$$\psi\circ(\xi\circ\mkern-2mu X) = \psi(\xi)\circ\mkern-2mu X.$$ 
This turns common stochastic operations into algebraic manipulations of predictable functions, which is both more effective and more compact than the standard calculus. It yields formulae such as $\e^{\zeta\sint X} = \Exp\bigs((\e^{\zeta\id}-1)\circ\mkern-2mu X\bigs)$ or $\Exp^\lambda(\xi\circ\mkern-2mu X) = \Exp(((1+\xi)^\lambda-1)\circ\mkern-2mu X)$ for $\lambda\in\N$. Unlike the canonical decompositions that commonly appear in classical stochastic calculus, semimartingale representations are measure-invariant.

Semimartingale representations allow a systematic evaluation of the predictable compensator (drift, $B^{\xi\circ\mkern-1mu X}$) of a represented process in terms of the predictable characteristics of the representing process $X$; this follows naturally from the \'Emery formula \eqref{eq:emery}. The drift calculation is further simplified by uniquely decomposing $X-X_0$ into two components, $X^\qc$ and $X^\ddp$, where $X^\qc$ is quasi-left-continuous and $X^\ddp$ equals the sum of its jumps at predictable times in the semimartingale topology. Only the drift of $\xi\circ\mkern-2mu X^\qc$ is evaluated via \eqref{eq:emery} since at predictable stopping times $\tau$ one has the simpler formula
$$\Delta B^{\xi\circ\mkern-1mu X}_\tau = \Delta B^{\xi\circ\mkern-1mu X^\ddp}_\tau = \E_{\tau_-}\bigsl[\xi_\tau(\Delta X_\tau)\bigsr].$$
In practice, $X^\qc$ is often an It\^o semimartingale. One may then rephrase the drift computation for the $\qc$ component in terms
of time rates, reverting to drift rates, quadratic variation rates (squared volatilities), and
jump intensities (L\'evy measures).

The rules of semimartingale representations together with drift evaluation give rise to the \emph{simplified calculus} of the title.  
In summary, the calculus provides a clear, systematic way to perform the ``$\circ$'' operations that we have showcased in the introduction and which we shall encounter again in various applications.

We shall now provide a rigorous setup of the paper. Below, we mostly rely on the notation of Jacod and Shiryaev \cite{js.03}. Throughout this section, let $m \in \N$ denote an integer.

\subsection{Preliminaries}\label{SS:2.1}
We explicitly shall allow quantities to be complex-valued. The reader interested only in real-valued calculus can easily always replace the general `$\Cx$--valued' by the special case `$\R$--valued' in their mind. We write  $\Cinf^m =\Cx^m \bigcup \{\NaN\}$ for some `non-number' $\NaN \notin \bigcup_{k\in\N}\Cx^k$. We introduce the functions  $\id:\Cinf^m\to\Cinf^m$  and $\hat \id: \Cinf^m \to \R^{2m} \cup \{\NaN\}$ by $\id(x)=x$  for all $x\in\Cinf^m$ and by
$$\hat \id(x) = (\Re x_1,\Im x_1,\ldots,\Re x_m,\Im x_m)^\top,\quad x \in \Cx^m;\qquad \hat \id(\NaN) = \NaN,$$
respectively. Observe that $\hat \id(x) \in \R^{2m}$ for $x \in \Cx^m$ contains the values of $\Re x$ and $\Im x$, interlaced. The introduction of $\NaN$ simplifies treatment of functions with undefined values, such as $\log(1+\id)$, that frequently arise in applications.

We fix a probability space $(\Omega ,\sigalg{F},\P)$ with a right-continuous filtration $\filt{F}$. We shall assume, without loss of generality, that all semimartingales are right-continuous, and have left limits almost surely. For a brief review of standard results without the assumption that the filtration is augmented by null sets, see Perkowski and Ruf \cite[Appendix~A]{Perkowski_Ruf_2014}. 

We denote the left-limit process of a (complex-valued) semimartingale $X$ by  $X_-$ and use the convention $X_{0-} = X_0$. We also set $\Delta X = X - X_-$.  We write $\hat{X} = \hat\id(X)$ and $X^\tau$ for $X$ stopped at some stopping time $\tau$.

For $\Cx$--valued semimartingales $X$ and $Y$ we set  
\begin{equation*}
 [X, Y]=[\Re X, \Re Y]-[\Im X, \Im Y] + i \left([\Re X,\Im Y]+[\Im X,\Re Y]\right).
\end{equation*}
If $X$ is $\Cx^m$--valued, then $[X, X]$ denotes the corresponding $C^{m \times m}$--valued quadratic variation, formally given by
\begin{equation*}
 [X, X]=  (I_m \otimes [1\ i]) \bigs[\hat X, \hat X\bigs] \left(I_m \otimes \left[ 
 \begin{array}{c}
1 \\ 
i%
\end{array}%
\right]\right) ,
\end{equation*}
where $I_m$ denotes the $m\times m$ identity matrix and $\otimes$ the Kronecker product. 
Furthermore, we write $[X, X]^c$ for the continuous part of the quadratic variation $[X, X]$.

As in \cite[II.1.4]{js.03}, we consider the notion of a predictable function on 
$\Oinf^m = \Omega \times [0,\infty) \times \Cinf^m$. Observe that every time-constant deterministic function, such as $\id$ or $\log(1+\id)$, can be considered a predictable function via the natural embedding of $\Cinf^m$ into $\Omega \times [0,\infty) \times \Cinf^m$. Let $\mu^X$ denote the jump measure of a semimartingale $X$ and $\nu^X$ its predictable compensator (under a fixed probability measure $\P$).  Then for a $\Cx$--valued  bounded predictable function $\xi$ with $\xi(0) = 0$ we  have
\[	\xi * \mu^X =  \xi  \Big(\hat \id^{-1}\Big) *  \mu^{\hat X} = \sum_{t \leq \cdot} \xi_t(\Delta X_t),
\]
 provided $|\xi| * \mu^X < \infty$.  If $X$ is special, we let  $B^X$ denote the predictable finite-variation part in the canonical decomposition of $X$, always assumed to start in zero, i.e., $B^X_0 = 0$. Recall from \cite[II.2.29a]{js.03} that
\begin{equation}\label{eq:190302.2}
\text{$X$ is special} \iff \left(|\id|^2 \wedge |\id|\right) * \nu^X < \infty.
\end{equation}

In this paper we call a time-constant deterministic function $h:\Cinf^m \rightarrow \Cinf^m$ a truncation function for $X$ if $|\id - h| * \mu^X < \infty$ 
and if 
\begin{equation}\label{eq:230315}
	X[h] = X - \left(\id - h \right) * \mu^X
\end{equation}
is special.

Next, let us  briefly discuss stochastic integrals.
Consider a $\Cinf^{1 \times m}$--valued process $\zeta$ and a $\Cx^m$--valued semimartingale $X$.  If $X$ is real-valued then we write $\zeta \in L(X)$ if both $\Re \zeta$ and $\Im \zeta$ are integrable with respect to $X$ (in the standard sense). We then write $\zeta \sint X = (\Re \zeta) \sint X + i (\Im \zeta) \sint X$. 
If $X$ is complex--valued,  then we say $\zeta \in L(X)$ if $( \zeta\otimes [1\ i]) \in L( \hat X)$, where $\otimes $ represents the Kronecker product. We then write 
\begin{equation*}
	\zeta \sint X = (\zeta\otimes [1\ i]) \sint \hat X
\end{equation*}
for the stochastic integral of $\zeta$ with respect to $X$. 

We sometimes shall work on stochastic intervals $\lc 0, \tau \lc$, where $\tau$ is a foretellable time, i.e., a stopping time that is almost surely equal to a predictable time.  An example is discussed in Subsection~\ref{SS:SExp}, where this setup allows to define the stochastic logarithm of a semimartingale that hits zero.  Let $\tau$ be a foretellable time with announcing sequence $(\tau_k)_{k \in \N}$, i.e., $\lim_{k \uparrow \infty} \tau_k = \tau$ and $\tau_k < \tau$ on $\{\tau > 0\}$.  Then we say a process $X$ is a semimartingale (local martingale, etc.) on $\lc 0, \tau \lc$ if $X^{\tau_k}$ is a semimartingale (local martingale, etc.) for each $k \in \N$.  We refer to Carr et al. \cite{CFR2011} and Larsson and Ruf \cite{larsson.ruf.20} for more details.

\subsection{Decomposition of a semimartingale into `continuous-time' and `discrete-time' components}
Denote by $\V$ the set of finite variation semimartingales; by $\V^{\d}$ the subset of $X\in\V$ such that 
$$X=X_0+\id * \mu^X;$$ 
and by $\V^{\d}_\sigma$ the set of semimartingales that belong sigma-locally to the class of pure-jump finite variation processes $\V^{\d}$. 

The following proposition recalls  a unique decomposition of a semimartingale $X$ into a semimartingale $X^\ddp$ that jumps at predictable times and a quasi-left-continuous semimartingale $X^\qc$. 
\begin{proposition}[\cite{cerny.ruf.21.bej}, Proposition~3.15]  \label{P:190729}
Every semimartingale $X$ has the unique decomposition 
\begin{equation*}
X = X_0 + X^{\qc} + X^{\ddp},
\end{equation*}
where  $X^\qc_0=X^{\ddp}_0 = 0$, $X^{\qc}$ is a quasi-left-continuous semimartingale,  $X^{\ddp}$ jumps only at predictable times, and $X^{\ddp} \in \V^{\d}_\sigma$. We then have $ [X^\qc,X^\ddp] = 0$.
\end{proposition}

\subsection{Stochastic exponentials and logarithms} \label{SS:SExp}
 If $X$ is a $\Cx$--valued semimartingale, then the stochastic exponential $\Exp(X)$ of $X$ is given by the formula (see \cite[Th\'eor\`eme~1]{doleans-dade.70})
\begin{equation}\label{eq:Exp}
 \Exp(X)=\e^{X-X_0-\frac{1}{2}[X,X]^c}\prod_{s\leq \cdot}\e^{-\Delta X_s}(1+\Delta X_s).
\end{equation}

In order to handle non-equivalent changes of measures, we extend the definition of stochastic logarithm (see \cite[II.8.3]{js.03}) to processes that can hit zero. 
To this end, for a $\Cx$--valued semimartingale $X$ define the stopping times%
\footnote{Since we have not assumed the filtration to be complete, the Debut theorem may not be applied. Hence, $\tau^X$ and $\tau^X_c$ themselves might not be stopping times. However, there always exists a stopping time almost surely equal to $\tau^X$ and a predictable time almost surely equal to $\tau^X_c$, respectively. Without loss of generality, we shall assume to work with such stopping times.}
\begin{align} 
\tau^X &= \inf \left\{t \geq 0: \inf_{s \leq t} |X_s| = 0\right\};\label{eq:221207.1}\\
\tau^X_c &=
\begin{cases}
\tau^X,&\qquad\lim_{t\uparrow \tau^X}X_t = 0 \\
\infty, &\qquad \lim_{t\uparrow \tau^X}X_t \neq 0
\end{cases}.  \label{eq:221207.2}
\end{align}
Here $\tau^X_c$ is the first time the running infimum of $\lvert X\rvert$ reaches zero continuously. Let us now additionally assume that $X$ is absorbed in zero if it ever hits zero. The stochastic logarithm $\Log(X)$ of $X$ is then given by
\begin{align*}
	\Log(X) = \frac{1}{X_-}\indicator{\{X_-\neq 0\}} \sint X \qquad \text{on }\lc 0,\tau^X_c\lc,
\end{align*}
where 
 $\frac{\indicator{\{X_{t-}\neq 0\}}}{X_{t-}}$ is defined to be zero on the set  $\{X_{t-}=0\}$, for all $t\geq 0$; see also Larsson and Ruf \cite{larsson.ruf.19}.

\subsection{Further details about predictable functions}
For this subsection, fix some $n \in \N$. For two predictable functions $\xi: \Oinf^m \rightarrow \Cinf^n$ and $\psi: \Oinf^n \rightarrow \Cinf$  we shall write $\psi(\xi)$ to denote the function $(\omega, t, x) \mapsto \psi(\omega, t, \xi(\omega, t, x))$ with the convention $\psi(\omega, t,  \NaN)=\NaN$. If $\psi$ and $\xi$ are predictable, then so is $\psi(\xi)$. 

For a predictable function $\xi: \Oinf^m \rightarrow \Cinf^n$ we shall write $\hat \xi = \hat \id(\xi)$ and $\xi^{(k)}$ for the $k$--th component of $\xi$, where $k \in \{1, \cdots, n\}$. We also write $\hat D \xi$ and $\hat D^2 \xi$ for the real derivatives of $\xi$, i.e., $\hat D_i \xi^{(k)}$ is the composition of the $i$--th element of the gradient of ${\xi^{(k)}}(\hat \id^{-1})$ and  the lift $\hat \id$ and $\hat D_{i,j}^2  \xi^{(k)}$ is the composition of the $(i,j)$--th element of the Hessian of $\xi^{(k)}( \hat \id^{-1})$ and the lift $\hat \id$, for $i,j \in \{1, \cdots, 2 m\}$.
Note that $\hat D \xi$ has dimension $n \times (2m)$, $\hat D^2 \xi$ has dimension $n \times (2m) \times (2m)$, and the domains of $\hat D \xi$, $\hat D^2 \xi$ equal that of $\xi$, i.e., $\Oinf^m$. If $\xi$ is analytic at a point, say $0$, then we also write $D \xi(0)$ and $D^2\xi(0)$ for the corresponding derivatives.

We want  to allow for predictable functions such as $\xi=\log(1+\id)$ 
whose effective domain is not the entire $\Cx$. To this end, we say that 
$$ \text{``a predictable function $\xi$ is \emph{compatible} with  $X$ if } \xi(\Delta X) \text{ is finite-valued}, \P\text{--almost surely.''}$$

\subsection{Semimartingale representation}\label{SS:2.5}
Often it will be useful to rely on representing a semimartingale with respect to another one. Such representations are worked out in 
\v{C}ern\'{y} and Ruf \cite{cerny.ruf.22.ejp}. Throughout this subsection let $X$ denote an $m$--dimensional semimartingale.

The following class $\Uni$ of predictable functions enjoys closedness with respect to common operations and a certain universality. A more general definition is possible (see \cite[Definition~3.2]{cerny.ruf.22.ejp}), but for our purposes, this universal class will suffice.
\begin{definition}[\cite{cerny.ruf.22.ejp}, Definition~3.4] \label{D:200518}
Let $\Uni^{n}$ denote the set of predictable functions $\xi: \Oinf^d \rightarrow \Cinf^n$ such that the following properties hold, $\P$--almost surely.
\begin{enumerate}[label={\rm(\arabic{*})}, ref={\rm(\arabic{*})}]
			\item\label{Uni:i}  $\xi_t(0) = 0$, for all $t \geq 0$.
			\item\label{Uni:ii}  $x \mapsto \xi_t(x)$ is twice real-differentiable at zero, 	for all $t \geq 0$.
			\item\label{Uni:iii} $\hat D \xi(0)$ and $\hat D^2 \xi(0)$ are locally bounded.
			\item\label{Uni:iv} There is a predictable locally bounded process $K>0$ such that 
	$$\sup_{K\abs{x}\leq 1} \frac{\mathopen{\bigs|} \xi(x)-\hat D \xi(0)\mkern2mu\hat{\id}(x)\mathclose{\bigs|}}{\abs{x}^2}\indicator{x\neq 0} \text{ is locally bounded.}$$
	\end{enumerate}
We write $\Uni=\bigcup_{n\in\N}\Uni^{n}$.
\qed
\end{definition}

\begin{definition}[\cite{cerny.ruf.22.ejp}, Definition~3.8]  \label{D:181101}
For a predictable function $\xi \in \Uni$ compatible with $X$ we use the notation 
	\begin{align*}
	\xi\circ\mkern-2mu X&{}=  \hat D \xi (0) \sint \hat X+ \frac{1}{2} \hat D^2 \xi(0) \sint \bigs[\hat X,\hat X\bigs]^c + \left(\xi -  \hat D \xi(0)\mkern2mu\hat\id  \right) * \mu^X. \tag*{\qed}
	\end{align*}	
\end{definition}
\noindent The following properties of semimartingale representations are worth pointing out.
\begin{itemize}
\item If $\xi\in\Uni$ is analytic or if $X$ is real-valued, then we may omit the hats on top of $D$, $D^2$, $X$, and $\id$ in the previous two definitions.

\item Using the notation of Proposition~\ref{P:190729}, we always have
\begin{equation}\label{eq:221103}
\Delta (\xi \circ\mkern-2mu X) =\xi(\Delta X); \qquad(\xi\circ X)^{\qc}=\xi\circ X^{\qc}; \qquad (\xi\circ X)^{\ddp}=\xi\circ X^{\ddp}.\end{equation} 
 
\item For sufficiently smooth $\xi:\R\to\R$, the partial sums 
$ \sum_{n\in\N}\xi \left(X^{t_n}-X^{t_{n-1}}\right)$ 
converge in ucp to $\xi\circ\mkern-2mu X$ as the time partition $(t_n)_{n\in\N}$ becomes finer; see \'Emery~\cite[Th\'eor\`eme~2a]{emery.78}.
\end{itemize}
\begin{example}[\cite{cerny.ruf.22.ejp}, Proposition~3.13(3)]\label{E:qv}
 We have $\id_i$, $\id_i\id_j  \in \Uni^1$, for all $i,j \in \{1,\ldots,m\}$, with
		\begin{align*}
			X^{(i)}     &= X^{(i)}_0+\id_i\circ\mkern-2mu X;\\
			\bigs[X^{(i)}, X^{(j)}\bigs] &= (\id_i\mkern2mu\id_j)\circ\mkern-2mu X. \tag*{\qed}
			\end{align*}
\end{example}

\begin{proposition}[\cite{cerny.ruf.22.ejp}, Proposition~3.13(6)] \label{P:200616}
 Let $Y$ be a predictable semimartingale of finite variation. Consider some $\xi \in \Uni$ compatible with $[X\,\, Y]^\top$ and assume $\xi(0, \cdot) = 0$.  
 Then $\xi(\cdot, \Delta Y)$ is in $\Uni$ and compatible with $X$. Furthermore,
	\begin{align*}
		\xi \circ (X, Y)=\xi(\,\cdot\, , \Delta Y) \circ X.
	\end{align*}
\end{proposition}
\begin{proposition}[Adapted from \cite{cerny.ruf.22.ejp}, Proposition~3.14]\label{P:integral}
Let $\zeta$ be a locally bounded $\mathbb{C}^{1\times m}$--valued predictable process. Then $\zeta\id \in \Uni$ and 
\[
\zeta\sint X=(\zeta\mkern2mu\id) \circ\mkern-2mu X.
\]
\end{proposition}

\begin{proposition}[\cite{cerny.ruf.22.ejp}, Proposition~3.15]\label{P:Ito}
Let $\U \subset \Cx^m$ be an open set such that $X_-, X \in \U$ and let $f: \U \rightarrow \Cx^n$ be  twice continuously real-differentiable.  
Then the predictable function $\xi^{f, X}:\Oinf^m\to\Cinf^n$ defined by
\begin{align*}
	 \xi^{f,X}(x) = 
		\begin{cases}
			f\left(X_- + x\right) - f\left(X_-\right), &\quad X_- + x \in \U\\
			\NaN, &\quad X_- + x \notin \U
		\end{cases},
		\qquad x\in\Cx^m,
\end{align*}
belongs to $\Uni^{n}$ and is compatible with $X$. Moreover,  
\begin{equation*}
f(X)=f(X_{0})+\xi^{f,X} \circ\mkern-2mu X. \tag*{\qed}
\end{equation*}
\end{proposition}

\begin{theorem}[\cite{cerny.ruf.22.ejp}, Theorem~3.17] \label{T:composition0}
The class $\Uni$ is closed under (dimensionally correct) composition, i.e., if $\xi \in\Uni^n$ and $\psi: \Oinf^n \rightarrow \Cinf$ is another predictable function with $\psi \in \Uni$ then $\psi(\xi)\in\Uni$. Furthermore, if $\xi\in\Uni$ is compatible with $X$ and $\psi\in\Uni$ is compatible with $\xi\circ\mkern-2mu X$, then $\psi(\xi)$ is compatible with $X$ and 
$$ \psi\circ(\xi\circ\mkern-2mu X) = \psi(\xi)\circ\mkern-2mu X.$$
\end{theorem}

\begin{proposition}[\cite{cerny.ruf.22.ejp}, Proposition~4.1]\label{P:190701b}
Assume $m\geq 2$. If  $\Delta X^{(2)}\neq -1$  then
\begin{align*}
\frac{\Exp \left(X^{(1)}\right)}{\Exp \left(X^{(2)}\right)} &= \Exp \left( \left(\frac{1+\id_1}{1+\id_2}-1\right) \circ\mkern-2mu X \right).
\end{align*}
Next assume instead that  $X^{(2)}_- \neq 0$, $X^{(2)} \neq 0$, $X^{(1)}$ does not go to zero continuously, and $X^{(1)}$ is absorbed in zero if it ever hits zero. Furthermore, let $\tau$ be the first time $X^{(1)} = 0$. Then
\begin{align*}
\Log \left(\frac{X^{(1)}}{X^{(2)}}\right) &=  \left(\frac{1+\id_1}{1+\id_2}-1\right)   \circ \bigs(\Log\bigs(X^{(1)}\bigs),\Log\bigs(X^{(2)}\bigs)^{\tau}\bigs).   
\end{align*}
\end{proposition}

\begin{proposition}[\cite{cerny.ruf.22.ejp}, Proposition~4.2]\label{P:190701}
We have
\begin{align} 
\Log \bigs(\e^X\bigs){}&{}=(\e^{\id}-1)\circ\mkern-2mu X; \label{eq:210410.1}\\
\lvert\Exp (X)\rvert{}&{}= \Exp\left(\left(\lvert 1+\id\rvert - 1\right) \circ\mkern-2mu X\right).\label{eq:190228.1}
\end{align}
If $\Delta X\neq -1$, then 
\begin{align} 
\Exp (X) &= \e^{\log (1+\id) \circ\mkern-1mu X};  \label{eq:220208}\\
\log\lvert\Exp (X)\rvert &= \log \lvert 1+\id\rvert \circ\mkern-2mu X,\notag
\end{align}
where $\log$ denotes the principal value logarithm.
\end{proposition}

\begin{proposition}[Adapted from \cite{cerny.ruf.22.ejp}, Proposition~4.3]\label{P:201126} 
Consider a $\Cx$--valued semimartingale $X$ and  $\alpha\in\Cx$. If $\alpha \in \Cx \setminus (0,\infty)$, assume furthermore $\Delta X \neq -1$.
We then have
\begin{equation}\label{eq:201126.1}
\mathopen{|}\Exp (X)\mathclose{|}^\alpha = \Exp \left( \left(\abs{1+\id}^\alpha -1\right)\circ\mkern-2mu X \right).
\end{equation}
\end{proposition}

Recall the notion of a truncated process $X[h]$ from \eqref{eq:230315}.
\begin{proposition}[\cite{cerny.ruf.22.ejp}, Proposition~5.6]\label{P:170822.1}
Fix $\xi\in\Uni$ compatible with $X$ and let $h$ (resp., $g$) be a truncation function for $X$ (resp., $\xi\circ\mkern-2mu X$). Then the following terms are well defined and the predictable compensator of $(\xi \circ\mkern-2mu X)[g]$ is given by 
\begin{align*}
B^{(\xi\circ\mkern-1mu X)[g]} &=\hat D \xi(0) \sint B^{\widehat{X[h]}}
		+ \frac{1}{2} \hat{D}^2 {\xi}(0) \sint \bigs[\hat X,\hat X\bigs]^c
		+ \left(g(\xi) -   \hat D \xi (0) \hat{h} \right) * \nu^X.   
	\end{align*}
If $\xi$ is analytic at $0$, $(\P \times A^X)$--almost everywhere, the following terms are well defined and
\begin{align*}
B^{(\xi\circ\mkern-1mu X)[g]} &=    D \xi(0) \sint B^{X[h]} 
		+ \frac{1}{2} D^2 {\xi}(0) \sint [X,X]^c
		+ \left(g (\xi) -   D \xi(0) h \right) * \nu^{X}.  
\end{align*}
\end{proposition}

\begin{remark}[\cite{cerny.ruf.22.ejp}, Remark~5.7]\label{R:two drifts}
Fix $\xi\in\Uni$ compatible with $X$. Recall that in the notation of Proposition~\ref{P:190729}, \eqref{eq:221103} gives
$$ \xi\circ X = \xi\circ X^\qc + \xi\circ X^\ddp.$$
Suppose now $\xi\circ\mkern-2mu X$ is special. One then has 
\begin{equation*}
	B^{\xi\circ\mkern-1mu X} =  B^{\xi\circ\mkern-1mu X^\qc} + B^{\xi\circ\mkern-1mu X^\ddp}.
\end{equation*}
Since the drift of $\xi\circ\mkern-2mu X^\ddp$ has a simple form given next, in practice Proposition~\ref{P:170822.1} is only used with $X=X^\qc$ to obtain $B^{\xi\circ\mkern-1mu X^\qc}$.
Indeed, the drift at predictable jumps times is given by
\begin{equation*}
 B^{\xi\circ\mkern-1mu X^\ddp} = \sum_{\tau\in\T_X}  \E_{\tau-}\bigsl[\xi_\tau(\Delta X_\tau)\bigsr]\indicator{\lc\tau,\infty \lc}.
\end{equation*} 
Here, $\T_X$  denotes a countable family of stopping times that exhausts the jumps of $X^\ddp$.
For each $X$, there are many ways to choose $\T_X$; it is sufficient to fix an arbitrary such family for each $X$.
 \qed
\end{remark}

\begin{corollary}[\cite{cerny.ruf.22.ejp}, Corollary~5.8]\label{C:170822.1}
Let $Y=Y_0+\xi\circ\mkern-2mu X$ for some $\xi \in \Uni$ compatible with $X$. Then the following holds.
\begin{gather*}
\begin{aligned}
	&\text{$\nu^Y$ is the push-forward measure of $\nu^X$ under $\xi$, that is, $\psi * \nu^Y = \psi(\xi) * \nu^X$}\\
	&\text{for all }\text{non-negative bounded predictable functions $\psi$ with $\psi(0)=0$.}
\end{aligned}
\end{gather*}
\end{corollary}

\begin{proposition}\label{P:190211} 
If $X$ has independent increments and if $\xi\in \Uni$ is compatible with $X$ and deterministic, then $\xi \circ\mkern-2mu X$, too, has independent increments. Moreover, if $X$ is a L\'evy process and if $\xi\in\Uni$ is compatible, deterministic, and time-constant, then $\xi \circ\mkern-2mu X$ is also a L\'evy process.
\end{proposition}
\begin{proof}
By \cite[II.4.15--19]{js.03}, $X$ has independent increments (respectively, is a L\'evy process) if and only if its characteristics are deterministic (respectively, the characteristics are absolutely continuous with respect to time and their Radon-Nikodym derivatives with respect to time are time-constant) relative to a truncation function for $X$. The claim now follows from Proposition~\ref{P:170822.1} and Corollary~\ref{C:170822.1}.
\end{proof}

\section{Multiplicative compensator} \label{S:3}
In many applications one seeks, for a given $\Cx$--valued semimartingale $Z$, a  predictable process 
of finite variation $K^Z$ that starts at $1$ and makes $\frac{Z}{K^Z}$ a local martingale. For example, changes of measure are frequently of the form $\frac{Z}{K^Z}$ where $Z$ is a real-valued process that after hitting zero remains in zero.  Another application arises when   $K^Z$ happens to be deterministic and $\frac{Z}{K^Z}$  to be a martingale. Then the multiplicative compensator is a device for computing expectations, namely,
$$\E\left[\frac{Z_t}{Z_0}\right]= K^Z_t, \qquad t \geq 0.$$
We will see in Section~\ref{S:4} that the L\'evy-Khintchin formula is but a special case of such setup.

Although we have in mind a situation where $Z$ has further structure, it transpires that one may express $K^Z$ directly in terms of $B^{\Log(Z)}$. This result (Theorem~\ref{T:mult comp}) is of independent interest because it simplifies and generalizes existing characterizations of multiplicative compensators; see  Jacod and Shiryaev \cite[II.8.21]{js.03}, Kallsen and Shiryaev \cite[Theorem~2.19]{kallsen.shiryaev.02},  and also L\'epingle and M\'emin \cite[Proposition~II.1]{lepingle.memin.78.ptrf}. Recall from \eqref{eq:221207.2} that $\tau^Z_c$ is the first time $Z$ reaches zero continuously.

\begin{theorem}[Multiplicative compensator]  \label{T:mult comp}
	Let $Z$ be a $\Cx$--valued semimartingale absorbed in zero if it ever hits zero. Assume that $\Log(Z)$ is special on $\lc 0, \tau_c^Z\lc$. 
	Assume next that 
	\begin{equation}  \label{eq:160816.4}
	\Delta B^{\Log(Z)} \neq -1 \qquad  \text{on  $\lc 0, \tau_c^Z \lc$.}  
	\end{equation}
Then the following statements hold.
\begin{enumerate}[label={\rm(\arabic{*})}, ref={\rm(\arabic{*})}]
\item\label{MC:2} We have $\Exp(B^{\Log(Z)})\neq 0$ on  $\lc 0, \tau_c^Z \lc$ and the process $M =\frac{Z}{\Exp(B^{\Log(Z)})} \indicator{\lc 0, \tau^Z_c\lc}$ is a local martingale on $\lc 0, \tau^Z_c\lc$. 
\item\label{MC:3} If $M$ is special (e.g., if $Z$ is special and $\liminf_{t \uparrow \tau_c^Z} |\Exp(B^{\Log(Z)})_t|  > 0$ 
	on $\{\tau_c^Z < \infty\}$),
	then $M$ is a local martingale on the whole positive real line. 
\end{enumerate}
\end{theorem}
Example~\ref{E:flight_to_cambodia} below illustrates how $M$ can fail to be a local martingale on the whole positive real line without the assumptions of \ref{MC:3}.

\begin{proof}[Proof of Theorem~\ref{T:mult comp}]
From \eqref{eq:160816.4} one obtains $\Log(\Exp(B^{\Log(Z)}))=B^{\Log(Z)}$ on $\lc 0, \tau_c^Z \lc$. Proposition~\ref{P:190701b} yields
\begin{equation}\label{eq:180703.2}
\begin{split}
		\Log(M) = \Log\left(\frac{Z}{\Exp\bigs(B^{\Log(Z)}\bigs)}\right) ={}& \frac{\id_1-\id_2}{1+\id_2}\circ \left(\Log(Z),B^{\Log(Z)}\right)   \qquad \text{on $\lc 0, \tau_c^Z \lc$}.
\end{split}
	\end{equation} 
Consequently, Proposition~\ref{P:170822.1} yields 
\begin{align*}
	B^{\Log(M)} &{}= B^{\Log(Z)} - B^{\Log(Z)} + \left(\frac{\id_1-\id_2}{1+\id_2} - \id_1 + \id_2\right) * \nu^{\left(\Log(Z), B^{\Log(Z)}\right) }\\
	&{}= \frac{\id_2^2 -\id_1 \id_2}{1+\id_2} * \nu^{\left(\Log(Z), B^{\Log(Z)}\right) }= \frac{\bigs(\Delta B^{\Log(Z)}\bigs)^2- \id \Delta B^{\Log(Z)}}{1+\Delta B^{\Log(Z)}} * \nu^{\Log(Z)} = 0 \qquad \text{on  $\lc 0, \tau_c^Z \lc$},
\end{align*}
which proves $M$ is a local martingale  on $\lc 0, \tau^Z_c\lc$. 

Assume now that $\liminf_{t \uparrow \tau_c^Z} |\Exp(B^{\Log(Z)})_t|  > 0$ on $\{\tau_c^Z < \infty\}$. As $\Exp(B^{\Log(Z)})\indicator{\lc 0, \tau_c^Z \lc}$  is predictable, we may assume by localization that $|\Exp(B^{\Log(Z)})\indicator{\lc 0, \tau_c^Z \lc}| > \delta$ for some $\delta > 0$; see \cite[Lemma~3.2]{larsson.ruf.20}. If additionally $Z$ is special, then clearly so is $M$.   
Let us now assume that $M$  is special. Then we may assume that $M$ is uniformly integrable.  Let now $(\tau_k)_{k \in \N}$ denote a non-decreasing sequence of stopping times such that $M^{\tau_k}$ is a uniformly integrable martingale and such that $\lim_{k \uparrow \infty} \tau_k= \tau_c^Z$.  With these localizations in place, we now fix $s,t \geq 0$ with $s < t$ and some $C \in \sigalg{F}_s$ and observe that
\[
	\E[M_t \indicator{C}] = \E\left[ \lim_{k \uparrow \infty} M_t^{\tau_k} \indicator{C}\right]
		=  \lim_{k \uparrow \infty}  \E\left[ M_t^{\tau_k} \indicator{C}\right] = \lim_{k \uparrow \infty}  \E\left[ M_s^{\tau_k} \indicator{C}\right] =  \E\left[ \lim_{k \uparrow \infty} M_s^{\tau_k} \indicator{C}\right] = \E[M_s \indicator{C}],
\]
proving the claim.
\end{proof}

\begin{remark}
The previous theorem shows that  $M=\frac{Z}{\Exp(B^{\Log(Z)})} \indicator{\lc 0, \tau^Z_c\lc}$ is a local martingale provided \eqref{eq:160816.4} holds, $Z$ is special, and  $\liminf_{t \uparrow \tau_c^Z} \Exp(B^{\Log(Z)})_t  > 0$. Conversely $Z = M \Exp(B^{\Log(Z)}) \indicator{\lc 0, \tau^Z_c\lc}$ yields the multiplicative decomposition of the semimartingale $Z$; see Jacod \cite{jacod.78} and M\'emin \cite{memin.78}.\qed
\end{remark}

When $Z$ is $\R$--valued, condition \eqref{eq:160816.4} makes sure that the expected percentage change in $Z$ is never equal to $-100\%$.  The next remark deals with the case where $\Re \frac{Z}{Z_-}$ is strictly positive. Note that $Z$ and $Z_-$ themselves may be $\Cx$--valued. 

\begin{remark}\label{R:170816.1}
 If  $Z_-\neq 0$ and $\Re \frac{Z}{Z_-}>0$, then \eqref{eq:160816.4} is automatically satisfied.
Indeed, for any predictable time $\tau$  we have 
	\begin{align*}
		0<\E_{\tau-} \left[\Re \frac{Z_\tau}{Z_{\tau-}}\right] = 1+ \Re \Delta B^{\Log(Z)}_\tau.
	\end{align*}
Hence,  $\Re \Delta B^{\Log(Z)}>-1$ and \eqref{eq:160816.4} holds. 
\qed
\end{remark}

\begin{remark}\label{R:181120.1}
When $Z$ can jump to zero, its multiplicative compensator on the interval $\lc \tau^Z, \tau^Z_c\lc$ is not defined uniquely. One may obtain another multiplicative compensator of $Z$ by replacing $\Exp(B^{\Log(Z)})$ in Theorem~\ref{T:mult comp}
 with another special semimartingale that is indistinguishable from $\Exp(B^{\Log(Z)})$ on the interval $\lc 0, \tau^Z\rc$ and satisfies a condition analogous to \eqref{eq:160816.4}.
This insight is used in the statement of Theorem~\ref{T:200603} and again in Corollary~\ref{C:PIIQ}.\qed
\end{remark}

 The next proposition contains some auxiliary results concerning the drift of the stochastic logarithm. These results contain sufficient conditions for the statements in Theorem~\ref{T:mult comp} to hold.
\begin{proposition}[Drift of stochastic logarithm] \label{P:181023}
Let $Z$ be a  $\Cx$--valued semimartingale absorbed in zero if it ever hits zero. If $Z$ is special then $\Log(Z)$ is special on $\lc 0, \tau_c^Z\lc$.  Moreover, if \eqref{eq:160816.4} holds then 
\[
	\left\{\liminf_{t \uparrow \tau_c^Z} \left| \Exp(B^{\Log(Z)})\right| > 0\right\}
		= \left\{ \liminf_{t \uparrow \tau_c^Z} \bigg(\Re B^{\Log(Z)^\qc}_t + \sum_{s \leq t} 
	\log\bigs\lvert 1+\Delta B^{\Log(Z)}_s\bigs\rvert\bigg) > -\infty\right\}.
\]
\end{proposition}

	As Example~\ref{E:flight_to_cambodia} below shows, 
	$\Log(Z)$ being special on $\lc 0, \tau_c^Z\lc$ in conjunction with \eqref{eq:160816.4} does not guarantee that $Z$ is special (on the whole positive real line).

\begin{proof}[Proof of Proposition~\ref{P:181023}]
	Observe that the process $\frac{1}{Z_-} \indicator{\{Z_- \neq 0\}}$ is locally bounded on  $\lc 0, \tau_c^Z\lc$. Hence, if $B^Z$ exists, then, with the help of Proposition~\ref{P:170822.1}, so does
\[
		 B^{\Log(Z)} = \frac{1}{Z_-} \indicator{\{Z_- \neq 0\}} \sint B^Z \qquad \text{on  $\lc 0, \tau_c^Z \lc $}.
\]
	This yields the first part of the statement.
	
	Assume now that \eqref{eq:160816.4} holds.  Then $\xi= \log |1+\id|$ is in $\Uni$  and compatible with $B^{\Log(Z)}$. We conclude from Proposition~\ref{P:190701} that  
	\begin{align*}
		\log \left\lvert\Exp \left(B^{\Log(Z)}\right)\right\rvert &= \log \lvert 1+\id\rvert \circ B^{\Log(Z)} = \Re B^{\Log(Z)} + (\log |1 + \id| -  \Re \id) * \mu^{B^{\Log(Z)}},
\end{align*}
	yielding the statement. 
\end{proof}

\begin{example}[Multiplicative compensator of an exponential process]  \label{E:190703.1}
	Assume that $Z$ is a special semimartingale of the form $Z =\e^{\xi \circ{} X}$ for some $\Cx^m$--valued semimartingale $X$ and compatible $\xi \in \Uni$.  Assume for simplicity that $\xi$ is analytic at $0$. Thanks to Proposition~\ref{P:190701} and Theorem~\ref{T:composition0} we have 
$\Log(Z) = ( \e^{\xi} - 1)\circ\mkern-2mu X$. Proposition~\ref{P:170822.1} now yields for $h= \id \indicator{\lvert\id\rvert  \leq 1}$ that 
\begin{equation}\label{eq:170903.1}
	\begin{split}
		B^{\Log(Z)} = B^{(\e^\xi-1)\circ\mkern-1mu X} ={}& D \xi(0) \sint B^{X[h]} + \frac{1}{2}  \left(D^2 \xi(0)  + D \xi(0)^\top D \xi(0) \right) \sint [X,X]^c \\
			&{}+ 
			 \left(\e^{\xi}  - 1 -  D \xi(0) h \right) * \nu^X.
	\end{split}
\end{equation}
This also gives 
$$
	\Delta B^{\Log(Z)} =  \int_{\Cx^m} \left( \e^{\xi(x)}  -  1 \right) \nu^X(\{\cdot\}, \d x)\neq -1. 
$$
Theorem~\ref{T:mult comp} now yields $\Exp(B^{(\e^\xi-1)\circ\mkern-1mu X})\neq 0$ and
\begin{equation*}
\frac{\e^{\xi\circ\mkern-1mu X}}{\Exp(B^{(\e^\xi-1)\circ\mkern-1mu X})}\text{ is a local martingale.}
\end{equation*}

Consider now the special case of the above with $\xi = \zeta\id$ for some locally bounded $\zeta$.%
\footnote{The following results also hold for general $\zeta \in L(X)$. However, for simplicity, in this paper we only discuss representations using the universal class $\Uni$, which requires the local boundedness of $\zeta$. For more details on generalizations, see the concluding Section~\ref{S:6}.} 
By Proposition~\ref{P:integral}, we have $ (\zeta\id) \circ\mkern-2mu X = \zeta \sint X$ while \eqref{eq:170903.1} simplifies to
	\begin{align*}
		 B^{\Log(Z)} & =B^{(\e^{\zeta\id}-1)\circ\mkern-1mu X}
		= \zeta \sint B^{X[h]} + \frac{1}{2} 
		\zeta^\top \zeta \sint [X,X]^c + 
			 \left( \e^{ \zeta\hspace{0.2mm} \id}  - 1  -  \zeta h \right) * \nu^X.
	\end{align*}	
Hence, by \eqref{eq:Exp} and Remark~\ref{R:two drifts} the multiplicative compensator of $Z$ equals
\begin{equation}\label{eq:190705.1}
\begin{split}
\Exp\Big(B^{\Log\left(\exp({\zeta\sint X})\right)}\Big) ={}&  \exp\left(\zeta \sint B^{X[h]^\qc} + \frac{1}{2} \zeta^\top \zeta \sint [X,X]^c  +  \left( \e^{\zeta\hspace{0.2mm}\id}  - 1  -  \zeta h \right) * \nu^{X^\qc}\right) \\
&{}\times \prod_{\tau\in\T_X}\E_{\tau-}\left[\e^{\zeta_\tau \Delta X_\tau}\right]\indicator{\lc\tau,\infty\lc}.
\end{split}
\end{equation}
Kallsen and Shiryaev \cite[Theorem~2.19]{kallsen.shiryaev.02} obtain \eqref{eq:190705.1} for real-valued $\zeta$ and $X$. In their work, 
the process $\log\bigs(\Exp\bigs(B^{\Log(\exp({\zeta\sint X}))}\bigs)\bigs)$ is called the exponential compensator of $\zeta \sint X$. \qed
\end{example}

\begin{example}[Multiplicative compensator of a power of a stochastic exponential]\label{E:201126}
Consider a $\Cx$--valued semimartingale $X$ and  $\alpha\in\Cx$. Assume $\Delta B^{\indicator{\id = -1} * \mu^X} \neq 1$. If $\alpha \in \Cx \setminus (0,\infty)$, assume furthermore $\Delta X \neq -1$. Finally, assume that $(|1+\id|^\alpha-1)\circ\mkern-2mu X$ is special. Then,
$$ \Exp(B^{(|1+\id|^\alpha-1)\circ\mkern-1mu X}) > 0$$
and
$$ \frac{|\Exp(X)|^\alpha}{\Exp(B^{(|1+\id|^\alpha-1)\circ\mkern-1mu X})} \text{ is a local martingale.}$$ 
This follows from Proposition~\ref{P:201126}, Theorem~\ref{T:mult comp}, and the composition rule in Theorem~\ref{T:composition0}.

A special case for $\alpha<0$ and a real-valued local martingale $X$ with $-1+\delta <\Delta X<1/\delta$ for some $\delta\in (0,1]$ appears in Kazi-Tani et al. \cite[Lemma~A.5]{kazi-tani.al.15}. A further special case, this time with $\alpha > 0$, appears in 
L\'epingle and M\'emin \cite[Proposition~II.3]{lepingle.memin.78.ptrf}.\qed
\end{example}

\begin{example}[$\Log(Z)$ special on $\lc 0, \tau^Z_c\lc$ but $Z$ not special on the whole time line]\label{E:flight_to_cambodia} 
Fix $\rho \in (0,  \e^{-1})$ and	consider the function
	\[
		f(t) = -\log(\rho - t)>1, \qquad 0\leq t < \rho.
	\]
	Here $\rho$ is the explosion time (to $\infty$) of $f$. Observe that $f$ satisfies $f'(t)=\e^{f(t)}$ for all $t \in [0, \rho)$.

Let now $V$ denote a  non-negative continuous local martingale with $\E[V_t] = 1$, $V_t \geq  \frac{1}{f(t)} $ for all $t \in [0, \rho)$, and  $V_t = 0$ for all $t \geq \rho$. Such a local martingale can be obtained, for example, by appropriately time-changing a Brownian motion started at one. Moreover, let $N$ denote an independent Poisson process with unit intensity.  Denote the stochastic process $(f(t))_{t\in [0,\rho)}$ by $f(\cdot)$ and define next the non-negative process
	\[
		Z = V \Exp\left(\e^{f(\cdot)}\indicator{\lc 0, \rho\lc}\sint N\right),
	\]
which is a semimartingale as the product of a continuous local martingale and a  process of finite variation.
	Then $\tau^Z_c = \tau^V_c =\rho$ and $\Log(Z)$ is special on $\lc 0, \tau^Z_c\lc$ with
	\[
		 B^{\Log(Z)} = B^{\e^{f(\cdot)} \sint N} = \int_0^\cdot \e^{f(t)} \d t = f(\cdot)-f(0) \qquad \text{on $\lc 0, \tau^Z_c\lc$}.
	\]
	Hence $\Exp(B^{\Log(Z)}) = \e^{f(\cdot) - f(0)}$ on $\lc 0, \tau^Z_c\lc$ and Theorem~\ref{T:mult comp} yields
	\[
		M = \e^{f(0) - f(\cdot)} {V} \Exp\left(\e^{f(\cdot)}\sint N\right)\indicator{\lc 0, \tau^Z_c\lc}   
	\]
is a local martingale on $\lc 0, \tau^Z_c\lc$ with jumps
	\[
		\Delta M = \e^{f(0)} {V}\Exp\left(\e^{f(\cdot)}\sint N\right)_-  \Delta N  \indicator{\lc 0, \tau^Z_c\lc}.
	\]
Hence, for each $t\geq 0$ one has 
$${\id^2} * \nu^M_t \leq {\id^2} * \nu^M_\rho  \leq \e^{2 f(0)} \prod_{0\leq t<\rho}\left(1+\e^{f(\cdot)}\Delta N_t\right)^2 \int_0^\rho {V}^2_t  \d t < \infty$$  and $M$ is special by \eqref{eq:190302.2}. However, 
 	\[
		\Delta Z  \, \geq \,   \e^{f(\cdot)} {V}  \Delta N  \indicator{\lc 0, \tau^Z_c\lc} \geq \frac{ \e^{f(\cdot)}}{f(\cdot) }  \Delta N  \indicator{\lc 0, \tau^Z_c\lc}
	\]
yields
\begin{align} \label{eq:190306}
	(\id^2 \wedge |\id|) * \nu^Z_\rho = {\id} * \nu^Z_\rho \geq \int_0^\rho \frac{\e^{f(t)}}{f(t)} \d t = \int_0^\rho \frac{f'(t)}{f(t)} \d t = \int_{f(0)}^\infty \frac{1}{u} \d u = \infty;
\end{align}
thus $Z$ is not special, again by \eqref{eq:190302.2}.

This now shows that requiring $M$ to be special is strictly weaker than requiring $Z$ to be special. This also gives an example 
where $\Log(Z)$ is special on $\lc 0, \tau^Z_c\lc$ but $Z$ is not special (on the whole positive line).

We now modify this example slightly to illustrate that $M$ in Theorem~\ref{T:mult comp} 
is not always a local martingale (on the whole positive line) if it is not required to be special \emph{a priori}. To see this, 
take $N$ to be a compound Poisson process with jumps of size $\pm 1$ with equal intensity, still independent of ${V}$. 
Now $B^{\Log(Z)} = 0$ on $\lc 0, \tau^Z_c\lc$, hence $M= Z$. However, as in \eqref{eq:190306}, 
the jumps of $M= Z$ are not locally integrable, hence $M$ cannot be a local martingale.
\qed
\end{example}

\section{Compensators of processes with independent increments}\label{S:4}

The following generalization of the L\'evy-Khintchin formula in Theorem~\ref{T:200603}\ref{T:200603.1} seems to be missing in the literature.  Kallsen and Muhle-Karbe \cite[Proposition~3.12]{kallsen.muhle-karbe.10.spa} prove a special case of Theorem~\ref{T:200603}\ref{T:200603.2}, assuming strict positivity (and real-valuedness) of the stochastic exponential, which allows an application of a measure change technique in their proof. In the same context, Cont and Tankov \cite[Proposition~8.23]{cont.tankov.04} do not require positivity but only treat the special case of L\'evy processes.

\begin{theorem}[Stochastic exponential of a process with independent increments] \label{T:200603}
Let $Y$ denote a $\Cx$--valued semimartingale with independent increments and define the deterministic time
\[
	\tau = \min \left\{t \geq 0: \P[\Delta Y_t = -1] = 1\right\}.
\]
Then  the following statements hold.
\begin{enumerate}[label={\rm(\arabic{*})}, ref={\rm(\arabic{*})}]
\item\label{T:200603.0} For $T\geq 0$ with $T < \tau$, $Y^T$ is special if and only if $\E[|\Exp(Y)_T|] < \infty$. Moreover, in this case we have 
\begin{align} \label{eq:flight:quatar}
	\E[\Exp(Y)_t] = \Exp(B^Y)_t ,\qquad 0\leq t\leq T.
\end{align}

\item\label{T:200603.1}  $Y$ is special on $\lc 0, \tau\lc$ if and only if $\E[|\Exp(Y)_t|] < \infty$ for all $t \geq 0$. Moreover, in this case we have 
$$\E[\Exp(Y)_t] = \Exp(B^Y)_t \indicator{\lc0, \tau\lc}(t),\qquad t \geq 0.$$
\item\label{T:200603.2} $\Exp(Y)$ is a local martingale if and only if it is a martingale.
\end{enumerate}
\end{theorem}

The proof of the theorem relies on the following lemma.

\begin{lemma}  \label{L:flight:Qatar}
Let $Y$ denote a $\Cx$--valued semimartingale with independent increments such that  $\P[\Delta Y_t = -1] < 1$ for all $t \geq 0$.  Then $\P[\Exp(Y)_t \neq 0] > 0$ for all $t \geq 0$.
\end{lemma}
\begin{proof}
Set  $U = -\indicator{\id=-1} \circ Y$. Then $U$ is special with  $\Delta B^U> -1$. Theorem~\ref{T:mult comp} and Remark~\ref{R:181120.1} now yield that $M = \frac{\Exp(U)}{\Exp(B^U)}$ is a local martingale. Since $M$ is bounded by the deterministic process $\frac{1}{\Exp(B^U)}$, we have $\E[M_t] = 1$ for all $t \geq 0$. This yields $\E[\Exp(U)_t]  = \E[M_t] \Exp(B^U)_t =\Exp(B^U)_t> 0$ for all $t \geq 0$. Since $\{\Exp(Y)_t \neq 0\} = \{\Exp(U)_t > 0\}$ for  all $t \geq 0$, this yields the statement.
\end{proof}

\begin{proof}[Proof of Theorem~\ref{T:200603}]
We first prove the assertion in \ref{T:200603.0}. To this end, fix $T \geq 0$ with $T < \tau$. Define next the process 
\begin{equation}\label{eq:200614.1}
	V = (|1+\id| - 1) \circ Y^T.
\end{equation}
By Proposition~\ref{P:201126}, we have $|\Exp(Y^T)\bigs|=\Exp(V)$ and by Proposition~\ref{P:190211}, $V$ has independent increments.  Moreover, by \eqref{eq:190302.2} and  Corollary~\ref{C:170822.1} we have that $Y$ is special if and only if $V$ is special.  Hence it suffices to argue the equivalence assertion with $Y^T$ replaced by $V$.

Assume now that $V$ is special. By assumption we also have  $\Delta B^V> -1$.
Thanks to Theorem~\ref{T:mult comp} and Remark~\ref{R:181120.1}, the process 
	$\frac{\Exp(V)}{\Exp(B^V)}$ is a  non-negative local martingale; in particular, its expectation is bounded by one.
	By \cite[II.4.15--19]{js.03}, $\Exp(B^{V})$ is deterministic,  we thus have 
	$\E[\Exp(V)_T] \leq \Exp(B^V)_T< \infty$, concluding the proof of the first implication. 
	
	Now assume that $\E[\Exp(V)_T] < \infty$. First, fix $t \in [0, T]$ and note that by Lemma~\ref{L:flight:Qatar} $\Exp(V)_T = \Exp(V)_t \Exp(V - V^t)_T$ is the product of two independent random variables, none of which is identically zero. Thus we get
	$$
	\E\left[\Exp(V)_T\right]  =  \E\left[\Exp(V)_t\right]\E\bigsl[\Exp(V -  V^t)_T\bigsr],
	$$
	which then yields $0<\E[\Exp(V)_t] < \infty$.
	Next, define the process 
	$$N =\left( \frac{\Exp(V)_u}{\E[\Exp(V)_u] }\right)_{u \geq 0}.$$ 
Fix now $s,t \geq 0$ with $s\leq t$ and $C \in \sigalg{F}_s$.  Then again by the independence of increments we have
	\[
		\E\left[N_t \indicator{C}\right] = \E\left[N_s \indicator{C}\right]  \E\left[\frac{\Exp(V -  V^s)_t}{\E\left[\Exp(V -  V^s)_t\right]}\right] = \E\left[N_s \indicator{C}\right],
	\]
	hence $N$ is a martingale. Denote by $\sigma$ the first time $V$ jumps by $-1$.
	The process $(W_u)_{u \geq 0}$, given by  
\[
	W_u = \E[\Exp(V)_u]  \indicator{\{u < \sigma\}} 
	+   \indicator{\{u \geq \sigma\}}  = \frac{\Exp(V)_u}{N_u}  \indicator{\{u < \sigma\}} +   \indicator{\{u \geq \sigma\}} 
\]
is a semimartingale uniformly bounded on compacts. In turn, this shows that $\Exp(V)=N W$ is special as the product of a martingale and a uniformly bounded on compacts semimartingale. Proposition~\ref{P:181023} now yields $V$ is special, proving the reverse implication in \ref{T:200603.0}. We also note that $\E[\Exp(V)_u] = \Exp(B^V)_u$ on paths where $u<\sigma$ by Theorem~\ref{T:mult comp}, hence for all $u\geq0$ by Lemma~\ref{L:flight:Qatar}.
 
Let us now consider the final assertion in \ref{T:200603.0}, namely \eqref{eq:flight:quatar}, provided $Y^T$ is special. We have already established that then also $0<\E[|\Exp(Y)_t|] < \infty$ for all $t \in [0, T]$.
Consider now the time $\rho = \inf \{t \in [0, T]: \Delta B^Y_t=-1\} $. Then $\rho$ is deterministic. 
Independence of increments yields  
	 $$\E[\Exp(Y)_t] = \E[\Exp(Y)_{\rho-}]\E[1+\Delta Y_{\rho}]\E[\Exp(Y-Y^\rho)_t] = 0 = \Exp(B^Y)_t,\qquad\qquad t \in [\rho,T];$$ 
	therefore we have $\E[\Exp(Y)_t] = \Exp(B^Y)_t$ for all $t \in [\rho, T]$.
Thus, without loss of generality, we may now just assume that $\rho = \infty$; in particular, that $\Delta B^Y \neq -1$ on $[0,T]$.  Then
Theorem~\ref{T:mult comp} and Remark~\ref{R:181120.1}
yield that $\frac{\Exp(Y)}{\Exp(B^Y)}$ is a local martingale. Its absolute value is bounded by
$$\frac{\Exp(V)}{\Exp(B^V)} \frac{\Exp(B^V)}{|\Exp(B^Y)|},$$ 
the product of a non-negative martingale and a deterministic semimartingale.  This shows that $\frac{\Exp(Y)}{\Exp(B^Y)}$ itself is a martingale, yielding the assertion.

For the statement in \ref{T:200603.1} note that $\Exp(Y) = 0$ on $\lc \tau, \infty\lc$. Also  $Y$ is special on $\lc 0,\tau\lc$ if and only if $Y^T$ is special for all $T\in [0,\tau)$. Then the statement follows from  the assertion in \ref{T:200603.0}.

	Finally, let us argue  \ref{T:200603.2} and let us assume that $\Exp(Y)$ is a local martingale. We may assume that $Y$ is constant after its first jump by $-1$, i.e., $Y=Y^\sigma$ with $\sigma = \inf\{t\geq 0:\Delta Y_t=-1\}$. Then $Y$ is a local martingale; hence $\E_{t-}[\Delta Y_t] = 0$ for all $t \geq 0$ by \cite[I.2.31]{js.03}, yielding $\tau = \infty$.  As above,  the local martingale  $\Exp(Y)$  is again bounded in absolute value by $\frac{\Exp(V)}{\Exp(B^V)} \Exp(B^V)$, the product of a non-negative martingale and a deterministic semimartingale.	This yields that  $\Exp(Y)$ is a martingale, concluding the proof.
\end{proof}

\begin{corollary}[L\'evy--Khintchin formula]\label{C:LKformula}
 Fix $u\in\R^{1\times d}$ and assume $X$ is an $\R^d$--valued semimartingale with independent increments starting at 0. Then
$$\E\bigsl[\e^{iuX_t}\bigsr]=\Exp\left(B^{\left(\e^{iu\id}-1\right)\circ\mkern-1mu X}\right)_t,\qquad t\geq 0.$$
Furthermore, if $h$ is a  truncation function for $X$ and if $X$ is a L\'evy process with 
drift rate $b^{X[h]}$ (relative to $h$) and jump measure $F^X$ one obtains 
$$\E\bigsl[\e^{iuX_t}\bigsr] =\exp\left(iub^{X[h]} t -\frac{1}{2} u\,[X,X]^c_t\, u\!^\top\! + t \int_{\R^d}(\e^{iux}-1-iuh(x))F^X(\d x) \right),\qquad t\geq 0.$$
\end{corollary}
\begin{proof}
Proposition~\ref{P:190701} yields $\e^{iuX} = \Exp\left(\left(\e^{iu\id}-1\right) \circ\mkern-2mu X\right)$. Since the jumps of $\e^{iuX}$ are bounded, an application of Theorem~\ref{T:200603}\ref{T:200603.1} yields the first statement. An application of Proposition~\ref{P:170822.1} then concludes.
\end{proof}

\begin{example}[Mellin transform of a signed stochastic exponential of an $\R$--valued process with independent increments] \label{E:190218}
	Fix $\alpha \in \Cx$ and for $j\in \{1,2\}$ let $f_j, \xi_j: \R\to\Cx$  denote the functions $ \xi_j = f_j(1+\id) - 1$  with
$$ f_1 = \lvert \id \rvert^\alpha\indicator{\id\neq 0};\qquad  
f_2 = \lvert \id \rvert^\alpha\left(\indicator{\id> 0}-\indicator{\id< 0}\right).$$
The functions $ \xi_1$ and $\xi_2$ are now extended to $\Oinf$ in the natural way by considering them to be constant in $t$, $\omega$, and in the imaginary component.
Note that $ \xi_1, \xi_2 \in \Uni^1$. Moreover, from Propositions~\ref{P:integral} and \ref{P:Ito}, and Theorem~\ref{T:composition0} we obtain, for $j\in \{1,2\}$ and any $\R$--valued semimartingale $Y$,
\begin{equation}\label{eq:200614.2}
\begin{split}
	\Log(f_j(\Exp(Y))) &= \frac{\indicator{\{\Exp(Y)_- \neq 0\}}}{ f_j(\Exp(Y)_-) }  \sint f_j(\Exp(Y))\\
		&= \frac{\indicator{\{\Exp(Y)_- \neq 0\}}}{ f_j(\Exp(Y)_-) }  \sint \left(\left( f_j\left(\Exp(Y)_-(1+\id)\right)  - f_j\left(\Exp(Y)_-\right) \right) \circ Y\right)\\
		&=  \xi_j  \indicator{\{\Exp(Y)_- \neq 0\}} \circ Y,
\end{split}
\end{equation}
yielding 
\begin{equation*}
f_j(\Exp(Y))=\Exp(\xi_j \circ Y),
\end{equation*}
which is a further generalization of \eqref{eq:190228.1} for real-valued $X$.%
\footnote{In this setting, $\xi_j\circ Y$ is defined for any $\Cx$--valued semimartingale $Y$ but the value of $\xi_j\circ Y$ is insensitive to the imaginary part of $Y$ by construction, so $Y$ is real-valued for all practical purposes. Observe that one can extend the functions $f_1$ and $f_2$ from $\R$ to $\Cx$ differently to get some action on the imaginary part of $Y$. For example, one could set, on $\Cx$,
$$ f_1 = \lvert \id \rvert^\alpha\indicator{\Re\id\neq 0};\qquad  
f_2 = \lvert \id \rvert^\alpha\left(\indicator{\Re\id> 0}-\indicator{\Re\id< 0}\right),$$
 keep the same definition of $\xi_1$ and $\xi_2$, and then extend to $\Oinf$ by making these functions constant in $t$ and $\omega$. 
In such case, the first two equalities in~\eqref{eq:200614.2} remain valid for arbitrary $\Cx$--valued $Y$ and the left-hand side  is  sensitive to $\Im Y$, but the third equality in~\eqref{eq:200614.2} still only works for real-valued $Y$. Hence, in the context of this example, it is reasonable to proceed with computations in \eqref{eq:191124.1} that are specialized to real-valued $Y$ and no longer hold for arbitrary $\Cx$--valued $Y$.}

Next, assume that $\xi_j \circ Y$ is special for $j\in \{1,2\}$; for example, this holds when $\Re\alpha = 0$ because the jumps of $\xi_j\circ Y$ are then bounded. Then, for $j\in \{1,2\}$, Proposition~\ref{P:170822.1} yields 
\begin{equation}\label{eq:191124.1}
\begin{split} 
B^{\xi_j \circ Y} ={}&   \alpha B^{Y[h]} +\frac{\alpha}{2}(\alpha-1)[Y,Y]^c
+ \left(\xi_j-\alpha h\right)*\nu^Y,
\end{split}
\end{equation}
where we may take $h=\id \indicator{|\id|\leq 1}$.
Assume now that $Y$ has independent increments. An application of Theorem~\ref{T:200603}\ref{T:200603.1} together with \eqref{eq:Exp} and Remark~\ref{R:two drifts} yields for $j \in \{1,2\}$ that
	\begin{align} \label{eq:200621}
		\E\left[f_j\left(\Exp(Y)_t\right)\right] = \Exp\bigs(B^{\xi_j\circ Y}\bigs)_t
		=\exp\Big(B_t^{\xi_j\circ Y^\qc}\Big)\prod_{s\leq t}\E\left[f_j(1+\Delta Y_s)\right], \qquad t\geq 0,
	\end{align}
	where from \eqref{eq:191124.1} one has
	\[
	B^{\xi_j\circ Y^\qc} = \alpha B^{Y[h]^\qc}+\frac{1}{2}\alpha(\alpha-1)[Y, Y]^c 	+(\xi_j-\alpha h)*\nu^{Y^\qc}.
	\]

From now on we shall fix $t \geq 0$ and acknowledge the explicit dependence on $\alpha$ by writing $f_j(\,\cdot\,;\alpha)$ and $\xi_j(\,\cdot\,;\alpha)$ for $j \in \{1,2\}$.
Next  define 
$$g_+(\alpha)=\E\left[|\Exp(Y)_t|^\alpha\indicator{\{\Exp(Y)_t>0\}}\right];\qquad g_-(\alpha)=\E\left[|\Exp(Y)_t|^\alpha\indicator{\{\Exp(Y)_t<0\}}\right].$$ 
From \eqref{eq:200621} we then have 
\begin{equation}\label{gplus}
2g_+(\alpha)=\E\left[f_1\left(\Exp(Y)_t;\alpha\right)\right] +  \E\left[f_2\left(\Exp(Y)_t;\alpha\right)\right]
 =\Exp\bigs(B^{\xi_1(\id;\alpha)\circ Y}\bigs)_t +  \Exp\bigs(B^{\xi_2(\id;\alpha)\circ Y}\bigs)_t,
\end{equation}
and similarly,
\begin{equation}\label{gminus}
2g_-(\alpha )= \E\left[f_1\left(\Exp(Y)_t;\alpha\right)\right] - \E\left[f_2\left(\Exp(Y)_t;\alpha\right)\right]
 =  \Exp\bigs(B^{\xi_1(\id;\alpha)\circ Y}\bigs)_t  - \Exp\bigs(B^{\xi_2(\id;\alpha)\circ Y}\bigs)_t.
\end{equation}
Next, define the following two conditional expectations:
\begin{align*}
	\text{if $\P[\Exp(Y)_t> 0] > 0$}: \qquad &\phi_+(u) = \E\left[\left. |\Exp(Y)_t|^{iu} \right| \Exp(Y)_t> 0\right], \qquad u \in \R;\\
	\text{if $\P[\Exp(Y)_t<0] > 0$} : \qquad &\phi_-(u)  = \E\left[\left. |\Exp(Y)_t|^{iu} \right| \Exp(Y)_t< 0\right], \qquad u \in \R.
\end{align*}
We can then compute
\begin{align*}
	\phi_+(u) = \frac{g_+(iu)}{g_+(0)}; \qquad \phi_-(u) = \frac{g_-(iu)}{g_-(0)}, \qquad u \in \R,
\end{align*}
provided $g_+(0)=\P[\Exp(Y)_t> 0] > 0$ and $g_-(0)=\P[\Exp(Y)_t< 0] > 0$. We have now obtained the Fourier transform of the the random variable $\log\lvert\Exp(Y)_t\rvert$ conditional on $\Exp(Y)_t\gtrless 0$. One is thus able to characterize the distribution of $\Exp(Y)_t$ via Mellin/Fourier inversion methods; see for example Galambos and Simonelli \cite{galambos.simonelli.04}.

Observe also that although the distribution of  
$\left\vert \Exp(Y)_t\right\vert$ conditional on $\Exp(Y)_t \gtrless 0$ corresponds to a strictly positive random variable, it cannot be thought of as a natural exponential of a process with independent increments (except in the trivial case when $\Exp(Y)_t> 0$ or $\Exp(Y)_t< 0$, $\P$--almost surely). Hence the characteristic functions $\phi_+$ and $\phi_-$ cannot be obtained from the classical L\'evy-Khintchin formula or from its generalization for processes with independent increments in \cite[Proposition~3.12]{kallsen.muhle-karbe.10.spa}. \qed
\end{example}

The next example illustrates the novelty of Example~\ref{E:190218} in a financial context. Throughout, $b^Z$ denotes the drift rate of a special L\'evy process $Z$.

\begin{example}\label{E:MV}
Let $X$ be a L\'evy process with characteristics $(b^{X[0]}=\mu,\sigma^2,\Pi= \lambda \Phi(0,\gamma^2))$, where 
$\Phi(\cdot,\cdot)$ denotes the cumulative normal distribution with a given mean and variance, respectively. Later we will use the specific numerical values  $\mu=0.2$, $\sigma=0.2$, $\lambda = 1$, and $\gamma = 0.1$, which are broadly consistent with the empirical distribution of the logarithmic returns of a well-performing stock. For simplicity we will use a zero risk-free rate.

Financial economics is concerned with optimal portfolio allocation over a period of time, e.g., $T=1$ year. Here we will consider optimality in the sense of mean--variance preferences. It is known that the optimal wealth process in this setting is given by $1-\Exp(-a(\e^{\id}-1)\circ\mkern-2mu X)$ where 
$$ a =\frac{b^{(\e^{\id}-1)\circ\mkern-1mu X}}{b^{(\e^{\id}-1)^2\circ\mkern-1mu X}} 
= \frac{\mu+\sigma^2/2+\lambda(\e^{\gamma^2/2}-1)}{\sigma^2+\lambda(\e^{2\gamma^2}-2\e^{\gamma^2/2}+1)}
\approx 4.48; $$
see, for example, Proposition~3.6, Lemma~3.7, Corollary~3.20, and Proposition~3.28 in \v{C}ern\'{y} and Kallsen \cite{cerny.kallsen.07}. Note that $a$ is the ratio of first and second moment of the arithmetic return of the stock.

\begin{figure}[t]
\centering%
\begin{subfigure}[b]{0.5\textwidth}
        		\centering
               \includegraphics[height=5.8cm]{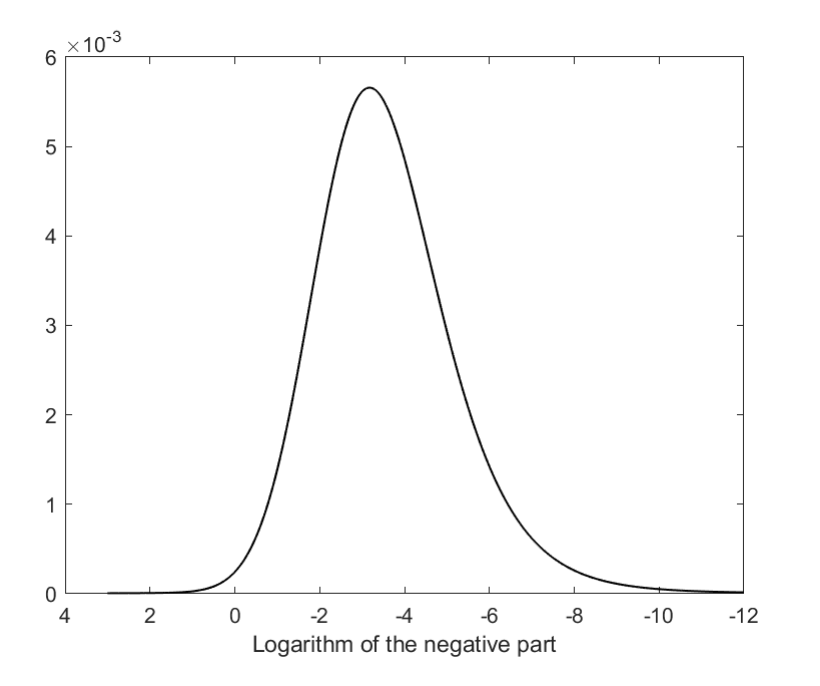}
                \caption{Subdensity of $\log \Exp(-a(\e^{\id}-1)\circ\mkern-2mu X)_T^-$.}
        \end{subfigure}%
\begin{subfigure}[b]{0.5\textwidth}
                \centering
                \includegraphics[height=5.8cm]{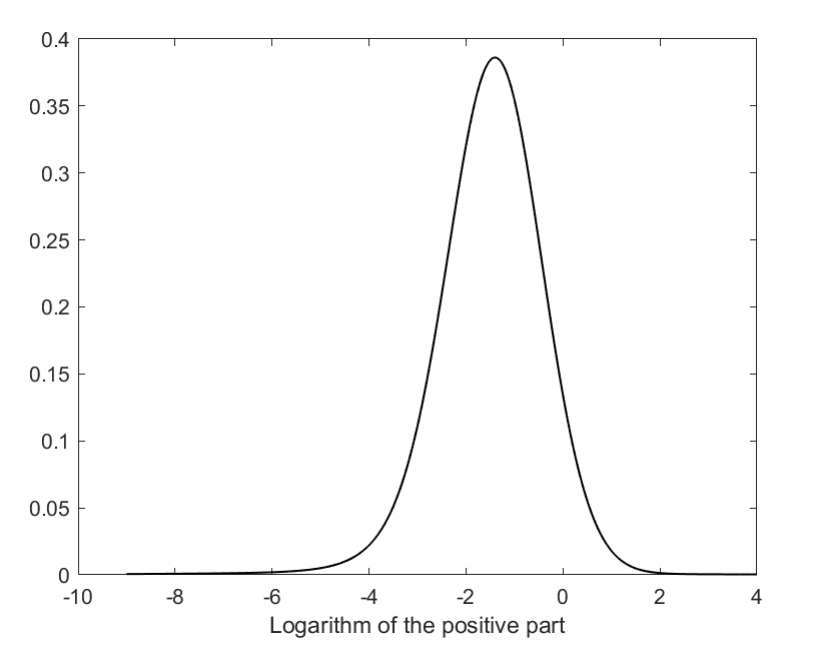}
                \caption{Subdensity of $\log \Exp(-a(\e^{\id}-1)\circ\mkern-2mu X)_T^+$.}
        \end{subfigure}
\caption{Distribution of a signed stochastic exponential}
\label{fig1}
\end{figure}

In practice, it is useful to know the distribution of the optimal terminal wealth 
$$1-\Exp(-a(\e^{\id}-1)\circ\mkern-2mu X)_T.$$ 
If the stochastic exponential $\Exp(-a(\e^{\id}-1)\circ\mkern-2mu X)$ is strictly positive, which is true in the empirically less important case $-1\leq a\leq 0$, this can be done by applying the L\'evy-Khintchin formula to the L\'evy process 
$$\log \Exp(-a(\e^{\id}-1)\circ\mkern-2mu X) = \log (1-a(\e^{\id}-1))\circ\mkern-2mu X$$ 
to obtain the characteristic function of the logarithm. In the commonly encountered situation with $a>0$, the model under investigation (and indeed all named L\'evy models used in finance) leads to a signed  stochastic exponential. The L\'evy--Khintchin formula is thus of no help but Example~\ref{E:190218} offers a way out.

Keeping the definitions of $\xi_1$, $\xi_2$, $g_+$, and $g_-$ from Example~\ref{E:190218}, we start by evaluating 
\begin{align*}
b^{\xi_1(\id;\alpha)\circ Y} ={}& b^{\xi_1(-a(\e^{\id}-1);\alpha)\circ\mkern-1mu X}=I_1(\alpha);\\
b^{\xi_2(\id;\alpha)\circ Y} ={}& b^{\xi_2(-a(\e^{\id}-1);\alpha)\circ\mkern-1mu X}=I_1(\alpha)-2I_2(\alpha),
\end{align*}
where
\begin{align*}
Y ={}& -a(\e^{\id}-1)\circ\mkern-2mu X;\\
I_1(\alpha) ={}& -\alpha a \left(\mu+\frac{1}{2}(1+a)\sigma^2\right) +\frac{1}{2}\alpha^2(a\sigma)^2+\int_{\R}(|1-a(\e^{x}-1)|^\alpha\indicator{a(\e^{x}-1)\neq 1}-1)\Pi(\d x);\\
I_2(\alpha) ={}& \int_{\R}|1-a(\e^{x}-1)|^\alpha\indicator{a(\e^{x}-1)> 1}\Pi(\d x).
\end{align*}
Next we obtain from \eqref{gplus} and \eqref{gminus}
\begin{align*}
g_+(\alpha) ={} \e^{I_1(\alpha)T}\frac{1+\e^{-2I_2(\alpha)T}}{2};\qquad
g_-(\alpha) ={} \e^{I_1(\alpha)T}\frac{1-\e^{-2I_2(\alpha)T}}{2}.
\end{align*}
Observe that without fixed jump times one has $g_+(0)>0$ and with $a>0$ also $g_-(0)>0$. In our setting we obtain $g_-(0)+g_+(0)=1$ with $g_-(0) \approx 2.2\%$ representing the probability that the stochastic exponential is negative in $T=1$ year. The conditional characteristic functions are clearly integrable, hence the standard density inversion formula can be applied. Figure~\ref{fig1} illustrates the subdensities of the logarithm of the negative and the positive part of the signed stochastic exponential $\Exp(-a(\e^{\id}-1)\circ\mkern-2mu X)_T$. Figure~\ref{fig2} shows the resulting distribution of the terminal wealth $1-\Exp(-a(\e^{\id}-1)\circ\mkern-2mu X)_T$ on the one-year horizon. \qed

\begin{figure}[t]
\centering%
\includegraphics[height=6.2cm]{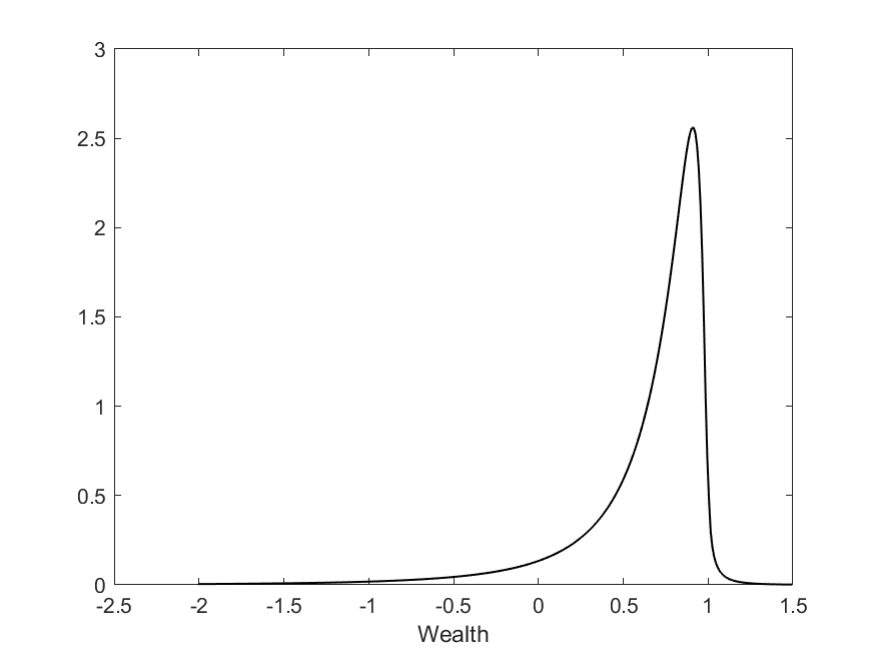}              
\caption{Density of the terminal wealth distribution $1-\Exp(-a(\e^{\id}-1)\circ\mkern-2mu X)_T$.}
\label{fig2}
\end{figure}

\end{example}

The next example illustrates how, in a general semimartingale model, the drift computation can be performed separately on predictable jump times.
\begin{example}[Multiplicative compensator calculation with predictable times of jumps]\label{E:181025.2}
Let $V$ denote a compound Poisson process with rate $\theta$ whose jumps have cumulative distribution function $F$.
Denote the jump times of $V$ by $(\rho_k)_{k \in \N}$ and set $\rho_0 = 0$. Let $L$ denote an independent special L\'evy process with drift rate $\mu$, variance rate $\sigma^2$, and jump measure $\Pi$ and set $X = L + V$. Next, let $\filt{F}$ denote the smallest right-continuous filtration such that $X$ is adapted and $\rho_k$ is $\sigalg{F}_{\rho_{k-1}}$--measurable, for each $k \in \N$.  Then $X$ is an $\filt{F}$--semimartingale and $\rho_k$ is $\filt{F}$--predictable, for each $k \in \N$.  Moreover,  its `continuous-time component' is precisely the L\'evy process; i.e., $X^\qc=L$. The `discrete times' lie in the set of predictable times of jumps, i.e., $\T_X=\bigcup_{k\in\N}\{\rho_k\}$.

We now interpret $X$ as the logarithmic price of an asset. 
Maximization of exponential utility calls for the multiplicative compensator of the utility process $Z = \e^{-\lambda\sint R}$ where $R =  \Log(\e^X)$ is the cumulative yield of $1\$$ investment in asset with price $\e^X$. 
Due to the presence of jumps at predictable times, the optimal investment strategy $\lambda$ will no longer be a constant dollar amount at all times, instead we will have one constant amount, say $\lambda_L$, on the `continuous' times, i.e., outside $\T_X$, and a different constant amount, say $\lambda_V$, on the `discrete' time set $\T_X$. 

By Proposition~\ref{P:190701} and Theorem~\ref{T:composition0} we obtain
$$\Log (Z) = \left(\e^{-\lambda\hspace{0.2mm} \id}-1\right)\circ R = \left(\e^{-\lambda(\e^{\id}-1)}-1\right)\circ\mkern-2mu X.$$
Observe that $\Log (Z)^{\qc}=\left(\e^{-\lambda(\e^{\id}-1)}-1\right)\circ\mkern-2mu X^\qc$ is a L\'evy process by Proposition~\ref{P:190211}.
For its drift rate $b^{\Log (Z)^\qc}$ we obtain by Proposition~\ref{P:170822.1} that
\begin{equation}\label{eq:181029.1}
 b^{\Log (Z)^\qc}=  - \lambda_L \mu + \frac{\sigma^2 \lambda_L}{2} \left(\lambda_L - 1 \right)  +\int_{\R}\left(
\e^{-\lambda_L \left( \e^x-1\right) }-1+\lambda_L x\right) \Pi(\d x). 
\end{equation}
This yields by \eqref{eq:Exp} and Remark~\ref{R:two drifts} that
\begin{align*}
	\Exp\bigs(B^{\Log(Z)}\bigs)_t &= \e^{b^{\Log (Z)^\qc}\!t}\,
	\prod_{\tau\in\T_X}\E_{\tau-}\left[\e^{-\lambda_V \left( \e^{\Delta X_\tau}-1\right) }\right] \indicator{\{\tau\leq t\}}
	\\
	&= 
\e^{b^{\Log (Z)^\qc}\!t}\,
\left(\int_{\R}\e^{-\lambda_V \left( \e^{x}-1\right) }F(\d x)\right)^{N_t},\qquad t \geq 0,
\end{align*}
where $N_t$ counts all predictable times of jumps in the interval $[0,t]$. 

Let us now fix a time horizon $T > 0$. 
By conditioning on $N$ and applying Theorem~\ref{T:200603}\ref{T:200603.1} we observe that $\E[Z_T]=\E[\Exp(B^{\Log(Z)})_T]$. This yields an explicit expression for the expected utility
\begin{align*}
\E\bigsl[\e^{-\lambda\sint R_T}\bigsr]={}&\E[Z_T]=\E\bigsl[\Exp\bigs(B^{\Log(Z)}_T\bigs)\bigsr]\\
={}&\e^{b^{\Log (Z)^\qc}\!T}\,
\e^{\kappa^\theta\left(\log\left(\int_{\R}\e^{-\lambda_V \left( \e^{x}-1\right) }F(\d x)\right)T\right)}
=\e^{b^{\Log (Z)^\qc}\!T}\, 
\e^{\theta\left(T\int_{\R}\e^{-\lambda_V \left( \e^{x}-1\right) }F(\d x)-1\right)},
\end{align*} 
where the drift rate $b^{\Log (Z)^\qc}$ \!is given in \eqref{eq:181029.1} and $\kappa^\theta=\theta(\e^{\id}-1)$ is the cumulant function of a Poisson variable with parameter $\theta$.\qed
\end{example}

\begin{example}\label{E:210411}
	In the setting of Example~\ref{E:190703.1}, assume that $\xi$ is deterministic. Thanks to Theorem~\ref{T:200603}\ref{T:200603.1} and \eqref{eq:Exp} we get
	\begin{align*}
		\E\left[\e^{\xi \circ\mkern-1mu X_t}\right] ={}& \exp\left(B_t^{\Log(Z)^\qc}\right)
		\prod_{s\leq t}\E\bigsl[\e^{\xi_s(\Delta X_s)}\bigsr],\qquad t\geq 0,
	\end{align*}
where $B_t^{\Log(Z)^\qc}$ is obtained from \eqref{eq:170903.1}.	See also Jacod and Shiryaev \cite[II.4.26]{js.03} for the special case when $\xi$ is time-constant and zero in a neighborhood of zero and $X$ is real-valued.  \qed
\end{example}

\section{Change of measure and its representation} \label{S:5} 
We now discuss how to compute drifts after an absolutely continuous change of measure. 
Since the collection of null sets of the new measure, say $\Qu$, may be larger than that of $\P$, one is compelled to study $\Qu$--drifts of processes that \emph{a priori} are not $\P$--semimartingales. The next proposition addresses this issue by offering a specific way to `lift' $\Qu$--semimartingales back up to $\P$. Its proof is provided at the end of this section.

\begin{proposition}\label{P:181111}
Let $M$ be a real-valued, non-negative uniformly integrable $\P$--martingale with $M_0=1$ and define the probability measure $\Qu$ by $\frac{\d\Qu}{\d \P} =M_\infty$.  
In line with \eqref{eq:221207.1} and  \eqref{eq:221207.2}, define the stopping time
\[
\tau_J^M = 
\begin{cases}
\tau^M,&\qquad\lim_{t\uparrow \tau^M}M_t \neq 0 \\
\infty, &\qquad \lim_{t\uparrow \tau^M}M_t = 0
\end{cases}. 
\]
Let $V$ be a $\Qu$--semimartingale. Then $V_{\tau_J^M-} = \lim_{t \uparrow \tau_J^M} V_t$ exists on $\{\tau_J^M < \infty\}$, 
$\P$--almost surely, and 
\[
	V_{\uparrow}=V \indicator{\lc 0, \tau_J^M\lc} + V_{\tau_J^M-} \indicator{\lc \tau_J^M, \infty\lc} 
\] 
is a $\P$--semimartingale  on $\lc 0, \tau_c^M \lc$ and $\Qu$--indistinguishable from $V$.
\end{proposition}

\begin{remark}\label{R:221208}
In view of Proposition~\ref{P:181111}, it entails no loss of generality to assume that a given $\Qu$--semimartingale is also a $\P$--semimartingale, at least on the open interval $\lc 0, \tau^M_c\lc$, where $M$ is the uniformly integrable martingale denoting the change of measure. This observation is relevant for Proposition~\ref{P:180706}, Theorem~\ref{T:Girsanov}, and Corollaries~\ref{C:Zurich}, \ref{C:181005}, and \ref{C:PIIQ} below.
Example~\ref{Ex:221208}\ref{Ex:221208.1} yields an instance where a $\Qu$--semimartingale is explicitly lifted to a $\P$--semimartingale.
\qed
\end{remark}

We now proceed to formulate a relevant version of Girsanov's theorem. For a stopping time $\rho$ and a predictable time $\sigma$, we say that process $X$ is a semimartingale (resp., a local martingale; special) on $\lc 0, \rho\rc\cap \lc 0, \sigma\lc$ if $X^\rho$ is a semimartingale (resp., a local martingale; special) on $\lc 0, \sigma\lc$.

\begin{proposition}[Girsanov's theorem]\label{P:180706}
Let $M$ be a real-valued, non-negative uniformly integrable $\P$--martingale with $M_0=1$ and define the probability measure $\Qu$ by $\frac{\d\Qu}{\d \P} = M_\infty$. For a $\P$--semimartingale $V$, the following are equivalent.
\begin{enumerate}[label={\rm(\roman{*})}, ref={\rm(\roman{*})}]
\item\label{P:180706:i} $V$ is $\Qu$--special.
\item\label{P:180706:ii} $V M$ is $\P$--special on  $\lc 0, \tau_c^M\lc$.
\item\label{P:180706:iii} $V +[V,\Log(M)]$ is $\P$--special on  $\lc 0, \tau^M\rc\cap \lc 0, \tau_c^M\lc$.
\end{enumerate} 
If either condition holds then the compensators corresponding to \ref{P:180706:i} and \ref{P:180706:iii} satisfy
\begin{equation}  \label{eq:221205}
	B^V_\Qu =  B^{V+[V,\Log(M)]} \qquad \text{on $\lc 0, \tau^M\rc\cap \lc 0, \tau_c^M\lc$}.
\end{equation}
Furthermore, the following are equivalent.
\begin{enumerate}[label={\rm(\roman{*}')}, ref={\rm(\roman{*}')}]
\item\label{P:180706:i'} $V$ is a $\Qu$--local martingale.
\item\label{P:180706:ii'} $V M$ is a $\P$--local martingale on  $\lc 0, \tau_c^M\lc$.
\item\label{P:180706:iii'} $V +[V,\Log(M)]$ is a $\P$--local martingale on $\lc 0, \tau^M\rc\cap \lc 0, \tau_c^M\lc$.
\end{enumerate} 
\end{proposition}
\begin{proof}
We will first argue the equivalence of \ref{P:180706:i}  and \ref{P:180706:ii} (and \ref{P:180706:i'}  and \ref{P:180706:ii'}, respectively) and then the equivalence of \ref{P:180706:ii} and \ref{P:180706:iii} (and \ref{P:180706:ii'}  and \ref{P:180706:iii'}, respectively).

Assume first that \ref{P:180706:i} holds, i.e., $V-B^V_\Qu$ is a local $\Qu$--martingale. 
As in \cite[III.3.8]{js.03} we then have $(V-B^V_\Qu)M$
is a local $\P$--martingale on $\lc 0, \tau_c^M \lc$. Since $B^V_\Qu M$ is  $\P$--special on  $\lc 0, \tau_c^M\lc$, we get \ref{P:180706:ii}. Assume now that  \ref{P:180706:ii} holds,  i.e., $VM-B^{VM}$ is a local $\P$--martingale on $\lc 0, \tau_c^M \lc$. Then $\frac{1}{M} (VM-B^{VM} ) = V- \frac{B^{VM}}{M}$ is a local $\Qu$--martingale. Since $\frac{B^{VM}}{M}$ is $\Qu$--special (as the product of the locally bounded process $B^{VM}$ and the $\Qu$--local martingale $\frac{1}{M}$), so is $V$, which yields 
\ref{P:180706:i}. Note that the same arguments also yield the equivalence of \ref{P:180706:i'}  and \ref{P:180706:ii'}.

We now argue the equivalence of \ref{P:180706:ii} and \ref{P:180706:iii}. 
Since $\frac{1}{M_-} \indicator{\{M_->0\}}$ is locally bounded on $\lc 0, \tau_c^M \lc$, \ref{P:180706:ii} is equivalent to $(\frac{1}{M_-} \indicator{\{M_->0\}}) \sint (VM)$ being $\P$--special on $\lc 0, \tau_c^M \lc$. Using integration by parts we get
\begin{align*}
\left(\frac{1}{M_-} \indicator{\{M_->0\}}\right)\sint (VM)  &= \left(\frac{V_-}{M_-} \indicator{\{M_->0\}}\right) \sint M  + \indicator{\{M_->0\}} \sint V+[\indicator{\{M_->0\}} \sint V,\Log(M)]  \\
	&= \left(\frac{V_-}{M_-} \indicator{\{M_->0\}}\right) \sint M  +  V+[V,\Log(M)]  
	\qquad \text{ on $\lc 0, \tau^M\rc\cap \lc 0, \tau_c^M\lc$.} 
\end{align*}
As $(\frac{V_-}{M_-} \indicator{\{M_->0\}}) \sint M$ is a  local $\P$--martingale on $\lc 0, \tau_c^M \lc$, we obtain the equivalence of \ref{P:180706:ii} and \ref{P:180706:iii}, and also of \ref{P:180706:ii'}  and \ref{P:180706:iii'}, respectively.

Assume now that \ref{P:180706:i}--\ref{P:180706:iii} hold. Using \ref{P:180706:i'} and \ref{P:180706:iii'} with $V$ replaced by $V-B_\Qu^V$ and noting that $[B_\Qu^V, \Log(M)]$ is a $\P$--local martingale  on  $\lc 0, \tau_c^M\lc$, we now get \eqref{eq:221205}.
\end{proof}

\begin{theorem}[Drift  after a change of measure] \label{T:Girsanov}
Let $Z$ be a semimartingale with $Z_0 = 1$, $\Log(Z)$ special on  $\lc 0,\tau^Z_c\lc$,  and $\Delta B^{\Log(Z)} \neq -1$ on $\lc 0,\tau^Z_c\lc$.  Assume that
$M=\frac{Z}{\Exp(B^{\Log(Z)})} \indicator{\lc 0,\tau^Z_c\lc}$ is a real-valued, non-negative uniformly integrable martingale and define the probability measure $\Qu$ by $\frac{\d\Qu}{\d \P} =M_\infty$. For a $\P$--semimartingale $V$, the following are equivalent. 
\begin{enumerate}[label={\rm(\roman{*})}, ref={\rm(\roman{*})}]
\item\label{T:Girsanov:i} $V$ is $\Qu$--special.
\item\label{T:Girsanov:ii} $V\frac{Z}{\Exp(B^{\Log(Z)})}$ is $\P$--special on $\lc 0, \tau_c^Z\lc$.
\item\label{T:Girsanov:iii} $V +  [V , \Log(Z)]$ is $\P$--special on $\lc 0, \tau^Z\rc\cap \lc 0, \tau_c^Z\lc$.
\end{enumerate}
If one of these conditions holds, then one has
\begin{equation} \label{eq:200610} 
B^V_\Qu = \frac{1}{1+ \Delta B^{\Log(Z)}} \sint B^{V +  [V , \Log(Z)]}  \qquad \text{on $\lc 0, \tau^Z\rc\cap \lc 0, \tau_c^Z\lc$}.
\end{equation}
Furthermore, the following are equivalent.
\begin{enumerate}[label={\rm(\roman{*}')}, ref={\rm(\roman{*}')}]
\item\label{T:Girsanov:i'} $V$ is a $\Qu$--local martingale.
\item\label{T:Girsanov:ii'} $V\frac{Z}{\Exp(B^{\Log(Z)})}$ is a $\P$--local martingale on  $\lc 0, \tau_c^Z\lc$.
\item\label{T:Girsanov:iii'} $V + [V,\Log(Z)]$ is a $\P$--local martingale on  $\lc 0, \tau^Z\rc\cap \lc 0, \tau_c^Z\lc$.
\end{enumerate}
\end{theorem}

\begin{proof}
Example~\ref{E:qv}, Theorem~\ref{T:composition0},  \eqref{eq:180703.2}, Proposition~\ref{P:200616}, and Proposition~\ref{P:integral} give, in this order, 
\begin{align}
V+[V,\Log(M)] ={}&\id_1(1+\id_2)\circ \left(V, \Log(M) \right)= \frac{\id_1(1+\id_2)}{1+\id_3}\circ \left(V, \Log(Z) ,B^{\Log(Z)}\right) \notag\\
={}& \frac{\id_1(1+\id_2)}{1+\Delta B^{\Log(Z)}}\circ \left(V, \Log(Z)\right)
= \frac{1}{1+ \Delta B^{\Log(Z)}} \sint (V +  [V , \Log(Z)]) \quad \text{on  $\lc 0, \tau_c^Z\lc$}.\label{eq:200616.2}
\end{align}
By Proposition~\ref{P:180706}, \ref{T:Girsanov:i} is equivalent to \ref{T:Girsanov:ii}, which in turn is equivalent to  $V +  [V , \Log(M)]$ being $\P$--special on $\lc 0, \tau^Z\rc\cap \lc 0, \tau_c^Z\lc$. By \eqref{eq:200616.2} and Lemma~4.2 of Shiryaev and Cherny \cite{shiryaev.cherny.02}, the latter is equivalent to  \ref{T:Girsanov:iii}.  
Proposition~\ref{P:180706}, identity \eqref{eq:200616.2}, and Lemma~4.2 of Shiryaev and Cherny \cite{shiryaev.cherny.02} also yield  \eqref{eq:200610}. The equivalence of \ref{T:Girsanov:i'}--\ref{T:Girsanov:iii'} is established similarly.
\end{proof}

The following statement is helpful when $V$ and $\Log(Z)$ are represented in terms of some common process $X$ because it delivers the same result as Theorem~\ref{T:Girsanov} without requiring the joint $\P$--characteristics of $V$ and $\Log(Z)$ as an input.
\begin{corollary}\label{C:200616}
Consider $Z$ and $\Qu$ as in Theorem~\ref{T:Girsanov} and suppose there are $\xi,\psi\in\Uni$ compatible with a $\P$--semimartingale $X$ such that $V = V_0+ \xi\circ\mkern-2mu X$ on $\lc 0,\tau^Z_c\lc$ and $\Log(Z) = \psi\circ\mkern-2mu X$ on $\lc 0,\tau^Z_c\lc$.
Then we have 
$$
B^{V^\qc}_\Qu ={} B^{\xi(1+\psi)\circ\mkern-1mu X^\qc};\qquad
B^{V^\ddp}_\Qu = \sum_{\tau\in\T_X} 
\frac{\E_{\tau_-}[\xi(\Delta X_\tau)(1+\psi_\tau(\Delta X_\tau))]}{\E_{\tau_-}[1+\psi_\tau(\Delta X_\tau)]}\indicator{\lc\tau,\infty\lc} \quad \text{on $\lc 0, \tau^Z\rc\cap \lc 0, \tau_c^Z\lc$}.
$$
\end{corollary}

\begin{proof}
Theorems~\ref{T:Girsanov} and \ref{T:composition0} yield 
$$B^V_\Qu = \frac{1}{1+ \Delta B^{\Log(Z)}} \sint B^{V +  [V , \Log(Z)]}=\frac{1}{1+ \Delta B^{\psi\circ X}} \sint B^{\xi(1+\psi)\circ X}
\qquad\text{on $\lc 0, \tau^Z\rc\cap \lc 0, \tau_c^Z\lc$}.$$ 
The statement now follows from Remark~\ref{R:two drifts}.
\end{proof}

The next example illustrates the convenience of Theorem~\ref{T:Girsanov} in a financial context when evaluating characteristic functions under a new measure. Throughout, $b^Z$ denotes the drift rate of a special L\'evy process $Z$.
\begin{example}\label{E:210410}
Let $S>0, S_->0$ be such that $\Log(S)$ is a L\'evy process. It is known from, e.g., Bender and Niethammer \cite{bender.niethammer.08} that under suitable conditions on the characteristics of $\Log(S)$, the absolutely continuous local martingale measure for $S$ whose density has the smallest $L^q(\P)$ norm is obtained by setting
$$ \frac{\d \Qu_q}{\d \P} 
= \frac{\Exp\bigs(\bigs(((1+\lambda_q\,\id)^+)^{\frac{1}{q-1}}-1\bigs)\circ \Log(S)\bigs)_T}
	{\Exp\Big(	 B^{\bigs(((1+\lambda_q\id)^+)^{\frac{1}{q-1}}-1\bigs)\circ \Log(S)}\Big)_T},\qquad  q>1,$$
where $\lambda_q\in\R$ solves 
$$b^{\id((1+\lambda_q\id)^+)^{\frac{1}{q-1}}\circ \Log(S)}=0.$$

To price contingent claims on $S$ under $\Qu_q$, it is helpful to know the characteristic function of $\log S$ under this measure. By Proposition~\ref{P:201126}, Theorem~\ref{T:200603} applied under $\Qu_q$, Corollary~\ref{C:200616}, and Example~\ref{E:210411} one obtains
\begin{align*}
\E_{\Qu_q}\bigs[\e^{\alpha (\log S_T-\log S_0)}\bigs]={}&\E_{\Qu_q}\bigs[\Exp\bigs(((1+\id)^\alpha-1)\circ \Log(S)\bigs)_T\bigs]\\
={}&\exp\Big(b^{((1+\id)^\alpha-1) ((1+\lambda_q\,\id)^+)^{\frac{1}{q-1}}\circ{} \Log(S) }T\Big)\\
={}&\exp\Big(b^{(\e^{\alpha\id}-1) ((1+\lambda_q(\e^{\id}-1))^+)^{\frac{1}{q-1}}\circ{} \log S}T\Big),\qquad \Re\alpha = 0,
\end{align*}
where the last equality follows from \eqref{eq:210410.1} and the composition rule in Theorem~\ref{T:composition0}.\qed
\end{example}

For completeness, we discuss the remaining two $\Qu$--characteristics of $V$. The second characteristic remains trivially unchanged. The third characteristic can be obtained from the following corollary.  Note that \cite[III.3.17]{js.03} provides an alternative expression to  \eqref{eq:181109} below.

\begin{corollary}[Predictable compensator under the new measure]\label{C:Zurich}
Consider $Z$ and $\Qu$ as in Theorem~\ref{T:Girsanov}. For a $\P$--semimartingale ${V}$, we have
\begin{align}  \label{eq:181109}
		\nu_\Qu^{{V}} =  \int_\R \frac{1+z} {1+\Delta B^{\Log(Z)}}  \nu^{{V, \Log(Z)}}(\cdot, \d z) \qquad \text{on $\lc 0, \tau^Z\rc\cap \lc 0, \tau_c^Z\lc$}.
\end{align}
	\end{corollary}
	
\begin{proof}
Let $n$ denote the dimension of $V$; i.e., assume that $V$ is an $\Cx^n$--valued process. Consider the set $G = G_1 \times G_2$ with $G_1 \subset \lc 0, \tau^Z\rc\cap \lc 0, \tau_c^Z\lc$ predictable and $G_2$ a closed set in $\Cx^n$ not containing a neighbourhood of zero.
Then  \eqref{eq:200610}, Example~\ref{E:qv}, Theorem~\ref{T:composition0}, and Proposition~\ref{P:170822.1} yield 
\begin{align*}
	\nu^{V}_\Qu (G) &=  B_\Qu^{\indicator{G}\circ V} 
	=
	\frac{1}{1+ \Delta B^{\Log(Z)}} \sint B^{\indicator{G}\circ V +  [\indicator{G}\circ V , \Log(Z)]} = 
		\frac{1}{1+ \Delta B^{\Log(Z)}} \sint
		B^{\indicator{G}(\cdot, \id_1) (1 + \id_2) \circ (V, \Log(Z))} \\
		&= \frac{1}{1+ \Delta B^{\Log(Z)}} \sint\left(\indicator{G}(\cdot, \id_1)(1 + \id_2) * \nu^{V, \Log(Z)} \right)\\
	&=  \indicator{G}(\cdot, \id_1)  \frac{1+\id_2 } {1+\Delta B^{\Log(Z)}} * \nu^{V, \Log(Z)}  =  \indicator{G}* \left(\int_\R \frac{1+z} {1+\Delta B^{\Log(Z)}}  \nu^{{V, \Log(Z)}} (\cdot, \d z)\right)
\end{align*} 
on $\lc 0, \tau^Z\rc\cap \lc 0, \tau_c^Z\lc$. This proves the statement.
\end{proof}

\begin{remark}
	In the setup of  Corollary~\ref{C:Zurich}, assume $V = V_0 + \xi^V \circ \Log(Z)$ on $\lc 0, \tau_c^Z\lc$ for some $\xi^V \in \Uni$ compatible with $\Log(Z)$. 
	Then \eqref{eq:181109} can be written as follows.
	\[
	\text{$ \nu_\Qu^{{V}}$ is the push-forward measure, under $\xi^V$, of }
		\frac{1+\id} {1+\Delta B^{\Log(Z)}}  \nu^{\Log(Z)} \qquad \text{on $\lc 0, \tau^Z\rc\cap \lc 0, \tau_c^Z\lc$}.
	\]
	To see this, consider first the formula \eqref{eq:181109} with $V=\Log(Z)$ and then apply Corollary~\ref{C:170822.1} under the measure $\Qu$.
\qed
\end{remark}

\begin{corollary}[Multiplicative compensator after a change of measure]\label{C:181005}
Consider $Z$ and $\Qu$ as in Theorem~\ref{T:Girsanov}. Assume $W$ is a $\Cx$--valued $\P$--semimartingale and $\Qu$--almost surely absorbed in zero if it eve hits zero. Moreover, suppose that $\Log(W)$ is special under $\Qu$ on $\lc 0, \tau^W_c\lc$ and
\[
	\Delta B^{\Log(W)}_\Qu \neq -1, \qquad \text{on $\lc 0, \tau^W_c\lc$,} \qquad \text{$\Qu$--almost surely.} 
\]
 Then $\Log(W) + [\Log(W), \Log(Z)]$ is $\P$--special on $\lc 0, \tau^Z\rc\cap \lc 0,\tau^W_c\wedge\tau^Z_c\lc$ and the multiplicative compensator of $W$ under $\Qu$ is given by
 \[
 	\Exp\left(B^{\Log(W)}_\Qu\right) 
	= \frac{\Exp\bigs(B^{\Log(Z)} + B^{\Log(W) + [\Log(W), \Log(Z)]}\bigs)}{\Exp\bigs(B^{\Log(Z)}\bigs) }
	\qquad \text{on $\lc 0, \tau^Z\rc\cap \lc 0,\tau^W_c\wedge\tau^Z_c\lc$}. 
 \]
\end{corollary}
\begin{proof}
The first assertion follows from Theorem~\ref{T:Girsanov}. Theorem~\ref{T:mult comp}\ref{MC:2} yields that $\Exp\left(B^{\Log(W)}_\Qu\right)$ is the $\Qu$--multiplicative compensator of $W$ on $\lc 0,\tau^W_c\lc$. Thanks to \eqref{eq:200610}, Proposition~\ref{P:integral}, and Proposition~\ref{P:200616}, we have
\begin{align*}
	\Exp\left(B^{\Log(W)}_\Qu\right) 
		&{}= \Exp\left(\frac{\id}{1+ \Delta B^{\Log(Z)}} \circ B^{\Log(W) +  [\Log(W) , \Log(Z)]}\right) \\
		&{}=\Exp\left(\frac{\id_1}{1+ \id_2} \circ \left(B^{\Log(W)  +  [\Log(W) , \Log(Z)]}, B^{\Log(Z)} \right) \right)\\
		&{}= \Exp\left(\left(\frac{1+\id_2+\id_1}{1+\id_2}-1\right) \circ \left(B^{\Log(W)  +  [\Log(W) , \Log(Z)]}, B^{\Log(Z)} \right) \right) 
\end{align*}
\text{on $\lc 0, \tau^Z\rc\cap \lc 0,\tau^W_c\wedge\tau^Z_c\lc$}, which  yields the claim by the generalized Yor formula in Proposition~\ref{P:190701b}.
\end{proof}

\begin{corollary}[Independent increments and change of measure]\label{C:PIIQ}
Consider $Z$ and $\Qu$ as in Theorem~\ref{T:Girsanov}. Assume $Z = \Exp(Y)$ for some $\P$--semimartingale $Y$ (hence 
$\tau^Z_c = \infty$). Consider a $\Cx$--valued $\P$--semimartingale $V$ and write
\[
	\tau = \min \left\{t \geq 0: \Qu[\Delta V_t = -1] = 1\right\}.
\]
Assume that $V$ and $Y$ have jointly independent increments under $\P$.   
Then $Y$ is $\P$--special, $V$ has independent increments under $\Qu$, and the following are equivalent for any time $T\in[0,\tau)$.
\begin{enumerate}[(i)]
\item\label{C:PIIQ.i} $\E_\Qu[|\Exp(V)_T|]<\infty$.
\item\label{C:PIIQ.ii} $V+[V,Y]$ is $\P$--special on $[0,T]$.
\end{enumerate}
Furthermore,  if one of these conditions holds, then
\[
	\E_{\Qu}[\Exp(V)_t] =  \frac{\Exp\bigs(B^{Y} + B^{V + [V, Y]}\bigs)_t}{\Exp\bigs(B^Y\bigs)_t },\qquad t\in[0,T].
\]
\end{corollary}
\begin{proof}
By Proposition~\ref{P:181023}, $Y^{\tau^Z}=\Log(Z)$ is $\P$--special. Therefore by \eqref{eq:190302.2},  
$$\left(|\id|^2 \wedge |\id|\right) * \nu^Y_u<\infty $$
on paths where $u\leq\tau^Z$, hence for all $u\geq0$ by Lemma~\ref{L:flight:Qatar} since $Y$ has independent increments. By \eqref{eq:190302.2} again, $Y$ is $\P$--special. 
 Thanks to Theorem~\ref{T:Girsanov} and Corollary~\ref{C:Zurich}, the $\Qu$--characteristics of $V$ are deterministic if 
	$V$ and $Y$ have jointly independent increments  under $\P$, hence $V$ indeed has independent increments under $\Qu$ by \cite[II.4.15--19]{js.03}.
	
By Theorem~\ref{T:200603}\ref{T:200603.0},	\ref{C:PIIQ.i} is equivalent to
\begin{enumerate}[label={\rm(\roman{*}')}, ref={\rm(\roman{*}')}]
\item\label{C:PIIQ.i'} $V$ is $\Qu$--special on $[0,T]$,
\end{enumerate}
which by Theorem~\ref{T:Girsanov} is equivalent to
\begin{enumerate}[resume*]
\item\label{C:PIIQ.ii'} $V + [V,\Log(Z)]$ is $\P$--special on $\lc 0 , \tau^Z\wedge T\rc$.
\end{enumerate}

Next, we shall establish equivalence of \ref{C:PIIQ.ii'} and \ref{C:PIIQ.ii}. Note that $V + [V,\Log(Z)] = V + [V,Y]$ on $\lc 0 , \tau^Z\wedge T\rc$. 
	 Hence, by \eqref{eq:190302.2},  
$$\left(|\id|^2 \wedge |\id|\right) * \nu^{V + [V,Y]}_u<\infty $$
on paths where $u\leq\tau^Z \wedge T$, hence for all $u\geq0$ by Lemma~\ref{L:flight:Qatar} since $V + [V,Y]$ has independent increments under $\P$.

Since $\Log(Z)=Y^{\tau^Z}$, we have $M=\frac{Z}{\Exp(B^{\Log(Z)})}=\frac{\Exp(Y)}{\Exp(B^Y)}$ by Remark~\ref{R:181120.1} and hence 
$$\E_{\Qu}[\Exp(V)_t]= \E\left[\Exp(V)_tM_t\right]=\E\left[\frac{\Exp(V)_t\Exp(Y)_t}{\Exp(B^Y)_t}\right]
=\E\left[\frac{\Exp(V+Y+[V,Y])_t}{\Exp(B^Y)_t}\right],\qquad t\in [0,T].$$  
Since both $Y$ and $V+[V,Y]$ are $\P$--special, Theorem~\ref{T:200603} completes the proof.
\end{proof}

\begin{example}\label{E:221123}
Theorem~\ref{T:intro3} in the introduction combines Theorem~\ref{T:Girsanov} and Corollary~\ref{C:181005} in a simplified setting. We shall now illustrate the usefulness of Theorem~\ref{T:intro3} on a calculation appearing in \cite{lepingle.memin.78.ptrf}.
Take $Y$ as in Theorem~\ref{T:intro3}. In addition, assume $Y$ is a local $\P$--martingale, i.e., $B^Y=0$, and $(1+\id)\log(1+\id)\circ Y$ is $\P$--special. The implication from \ref{T:intro3:iii} to \ref{T:intro3:i} in Theorem~\ref{T:intro3} with $X = \log(1+\id)\circ Y$ now yields that $\log(1+\id)\circ Y$ is $\Qu$--special with compensator $B^{\log(1+\id)\circ Y}_\Qu=B^{(1+\id)\log(1+\id)\circ Y}$. Next, fix some $\lambda\in(0,1)$ and observe that $\Exp^\lambda(Y)=\Exp(((1+\id)^\lambda-1)\circ Y)$ is $\P$--special, hence $\Exp^{\lambda-1}(Y)$ is $\Qu$--special by the implication from \ref{T:intro3:ii} to \ref{T:intro3:i} in Theorem~\ref{T:intro3}. 
Now  \eqref{eq:220208} yields
$$\e^{-(1-\lambda)(\log(1+\id)\circ Y-B^{(1+\id)\log(1+\id)\circ Y})} = \Exp^{\lambda-1}(Y)\e^{(1-\lambda)B^{(1+\id)\log(1+\id)\circ Y}}, $$
which by Jensen's inequality 
is a local $\Qu$--submartingale; hence its $\Qu$--multiplicative compensator is non-decreasing. By the implication from \ref{T:intro3:i'} to \ref{T:intro3:ii'} in Theorem~\ref{T:intro3}, this $\Qu$--multiplicative compensator coincides with the $\P$--multiplicative compensator of 
$$\Exp^{\lambda}(Y)\e^{(1-\lambda)B^{(1+\id)\log(1+\id)\circ Y}}, $$
hence by Theorem~\ref{T:intro2} the $\Qu$--multiplicative compensator equals
$$ \Exp(B^{((1+\id)^\lambda-1)\circ Y})\e^{(1-\lambda)B^{(1+\id)\log(1+\id)\circ Y}}.$$
The non-decreasing property deduced earlier now yields
$$ \Exp(B^{((1+\id)^\lambda-1)\circ Y})\geq \e^{(\lambda-1)B^{(1+\id)\log(1+\id)\circ Y}},$$
which is the statement of a key inequality in \cite[Lemma~III.4, Equation~(3.4)]{lepingle.memin.78.ptrf}.\qed
\end{example}

We conclude with an explicit example of drift calculation after a non-equivalent change of measure where $Z$ is allowed to attain zero continuously.

\begin{example}  \label{Ex:221208}
Let $Z$ and $\Qu$ be as in Theorem~\ref{T:Girsanov}.
 We shall compute the $\Qu$--drift of two $\Qu$--semimartingales, $V$ and $U$. \begin{enumerate}[label={\rm(\roman{*})}, ref={\rm(\roman{*})}]
\item \label{Ex:221208.1}	We first consider the $\Qu$--semimartingale $V = \Log(\frac{1}{Z})$. 
Proposition~\ref{P:190701b} yields the representation 
$$ V = \Log\left(\frac{1}{Z}\right)=\left(\frac{1}{1+\id}-1\right)\stackrel{\Qu}{\circ}\Log(Z),$$
where we use $\stackrel{\Qu}\circ$ to emphasize that this representation only holds under $\Qu$. 
The lifted version of $V$ from Proposition~\ref{P:181111} reads
\begin{equation*}
	 V_{\uparrow} =  \left(\frac{1}{1 + \id}-1 \right)\indicator{\id \neq -1  } \circ \Log(Z)\qquad\text{on  $\lc 0, \tau_c^Z \lc$.}
\end{equation*}
Then Corollary~\ref{C:200616} with $\xi = (\frac{1}{(1 + \id)}-1) \indicator{\id \neq -1  }$ and $\psi=\id$ yields
\begin{align*}
	B^{\Log(\frac{1}{Z})}_\Qu &{}= B^{(-\indicator{\id = -1  }-\id)\circ \Log(Z)^\qc}
	-\sum_{\tau \in \T_Z} \frac{\E_{\tau_-}[\indicator{\{\Delta \Log(Z)_\tau = -1\}}]+\Delta B^{\Log(Z)}_\tau} {1+\Delta B^{\Log(Z)}_\tau} 
	\indicator{\lc\tau,\infty\lc},
\end{align*}
$\Qu$--almost surely. In particular, if $Z$ is a uniformly integrable martingale, then $B^{\Log(Z)} = 0$ and
\begin{align*}
	B^{\Log(\frac{1}{Z})}_\Qu =  -\indicator{\id=-1} * \nu^{\Log(Z)}, \qquad \text{$\Qu$--almost surely}.
\end{align*}
\item Let us now consider a second example in this setup. Assume that   $\id^2 * \nu^Z< \infty$ and consider the process 
$$U = \frac{1}{2} \left(Z^2 - 1- [Z,Z] \right) =  \frac{1}{2}  \left((Z_-+\id)^2-Z_-^2-\id^2\right)\circ Z = (Z_-\id)\circ Z.$$
Note that $U$ is a $\P$--semimartingale. 
Corollary~\ref{C:200616} with $\xi = Z_-\id$ and $\psi=\frac{\id\indicator{\{Z_-\neq 0\}}}{Z_-}$ yields
\begin{equation*}
\ \quad B^{U}_\Qu = B^{(Z_-\id+\id^2)\circ Z^\qc} +
	\sum_{\tau \in \T_Z} \frac{\E_{\tau_-}[(\Delta Z_\tau)^2]+Z_{\tau-}\Delta B^{Z}_\tau} {1+\frac{\Delta B^Z_\tau}{Z_{\tau-}}} \indicator{\lc\tau,\infty\lc},\qquad
	\text{$\Qu$--almost surely,}
\end{equation*}
where we have used $\id^2\indicator{\{Z_-= 0\}}\circ Z=0$. In particular, if $Z$ is a uniformly integrable martingale, then $B^{Z} = 0$ and
\begin{equation*}
	B^{U}_\Qu =  B^{[Z,Z]}, \qquad \text{$\Qu$--almost surely}.\tag*{\qed}
\end{equation*}
\end{enumerate}
\end{example}
We conclude this section with a proof of Proposition~\ref{P:181111}.
\begin{proof}[Proof of Proposition~\ref{P:181111}] 
By localization, we may and shall assume, without loss of generality, that $\tau_c^M = \infty$. 

Consider the process
\[
	W = V \indicator{\lc 0, \tau_J^M\lc} + \left(\limsup_{t \uparrow \tau_J^M} V_t\right) \indicator{\lc \tau_J^M, \infty\lc}. 
\] 
Note that $W$ is a $\Qu$--semimartingale, $\Qu$--indistinguishable from $V$, since $\Qu[\tau_J^M = \infty]=1$. Define the two processes $\underli{W}=(\underli W_t)_{t \geq 0}$ and $\widebar{W}=({\widebar W}_t)_{t \geq 0}$ by $\underli W_0 = W_0 = {\widebar W}_0$, $\underli W_t = \liminf_{s \uparrow t} W_s$, and ${\widebar W}_t = \limsup_{s \uparrow t} W_s$ for all $t > 0$.  Then $\underli{W}$ and $\widebar{W}$ are predictable. Hence, the first time $\rho$ that $W$ fails to be left-continuous, namely the first time when $\underli{W}$ does not equal $\widebar{W}$ or is not real-valued, is a predictable time. 
Note that
\[
	\P[\rho<\tau^M] = \E\left[\indicator{\{\rho<\tau^M\}} \frac{1}{M_\rho} M_\rho\right]
		= \E_\Qu\left[\indicator{\{\rho<\tau^M\}}\frac{1}{M_\rho} \right] = \E_\Qu\left[\indicator{\{\rho<\infty\}}\frac{1}{M_\rho} \right] = 0
\]
since $W$ is a $\Qu$--semimartingale and hence $\Qu[\rho<\infty]=0$. Moreover, since $W$ is constant on the interval $\lc \tau_J^M, \infty\lc$, we have $\rho= \tau^M_J$ on $\{\rho<\infty\}$, $\P$--almost surely. 
Since $\rho$ is predictable, by \cite[I.2.31]{js.03} we have 
	\[
		0  = \E_{\rho-}[\Delta M_{\rho}] = \E_{\rho-}\left[\Delta M_{\rho} \indicator{\{\rho = \tau^M_J\}}\right] 	
		\qquad \text{on $\{\rho< \infty\}$.}
	\]  
Since $\Delta M_\rho < 0$ on $\{\rho = \tau^M_J\}$, this yields  $\P[\rho  = \tau^M_J]=0$; hence we have $\P[\rho = \infty] = 1$.
This shows that $W$ is $\P$--almost surely left-continuous; in particular, $\lim_{t \uparrow \tau_J^M} V_t$ exists on $\{\tau_J^M<\infty\}$ and $W=V_\uparrow$. 

We also note that $V_\uparrow$ has right-continuous paths. Indeed, we use $\rho$ again, but now to denote the first time that $V$ fails to be right-continuous. Then $\rho$ is a stopping time. As above, since $V$ is a $\Qu$--semimartingale we have $\rho \geq \tau^M$. Since $V_\uparrow$ is constant after time $\tau^M$, we indeed have $\rho = \infty$, yielding that $V_\uparrow$ has right-continuous paths with left limits, $\P$--almost surely.

The $\Qu$--semimartingale $V$ can be written as the sum of a $\Qu$--local martingale $V^{\rm{lm}}$ and a $\Qu$--semimartingale $V^{\rm{fv}}$ that is $\Qu$--almost surely of finite variation. To show the proposition, it now suffices to argue that their corresponding lifts   $V^{\rm{lm}}_\uparrow$ and $V^{\rm{fv}}_\uparrow$ are $\P$--semimartingales.

Let us first consider $V^{\rm{fv}}_\uparrow$.  Denote by $\rho$ the first time that $V^{\rm{fv}}_\uparrow$ is of infinite variation. Then $\rho$ is a predictable time. Using the same arguments as above we can argue that $\rho \geq \tau^M$, $\P$--almost surely, then that $\rho = \tau^M$ on $\{\rho < \infty\}$, $\P$--almost surely, and then  that indeed $\rho = \infty$, $\P$--almost surely.  Thus, $V^{\rm{fv}}_\uparrow$ is right-continuous with left limits, and $\P$--almost surely of finite variation, hence is a $\P$--semimartingale.

Finally, let us consider $V^{\rm{lm}}_\uparrow$, which is $\Qu$--indistinguishable from $V^{\rm{lm}}$.  Let $(\rho_n)_{n \in \N}$ denote a localization sequence so that $(V^{\rm{lm}}_\uparrow)^{\tau_n}$  is a $\Qu$--martingale and $\lim_{n \uparrow \infty} \rho_n = \infty$, $\Qu$--almost surely. Without loss of generality, we may assume that $\rho_n \leq n$, $\P$--almost surely. Next, let us define the stopping times 
	\[
		\rho_{n\uparrow} = \rho_n \indicator{\{ \rho_n < \tau^M\}} + (\rho_n + n) \indicator{\{ \rho_n \geq \tau^M\}}, \qquad n \in \N.
	\]
Then  $\Qu[\rho_{n\uparrow} = \rho_n] = 1$ for all $n \in \N$ since $\Qu[\tau^M = \infty] = 1$. Hence $(\rho_{n\uparrow})_{n \in \N}$ is again a localization sequence for the $\Qu$--local martingale $V^{\rm{lm}}_\uparrow$. Moreover,  the limit $\rho = \lim_{n \uparrow \infty} \rho_{n\uparrow}$ exists and satisfies $\rho =\tau^M$ on $\{\rho < \infty\}$, $\P$--almost surely; hence  $\P[\rho_{n\uparrow} < \rho]=1$ for each $n \in \N$. Therefore $\rho$ is $\P$--almost surely equal to a predictable time and as above we have again  $\P[\rho = \infty] = 1$. Now \cite[III.3.8c]{js.03} yields that $M V^{\rm{lm}}_\uparrow$ is a $\P$--local martingale. Dividing this process by the strictly positive $\P$--semimartingale $M + \indicator{\lc \tau^M, \infty\lc}$ and adding the finite-variation process $V_{\tau_J^M-}^{\rm{lm}} \indicator{\lc \tau^M, \infty\lc}$ yields that $V^{\rm{lm}}_\uparrow$ is a $\P$--semimartingale as claimed.
\end{proof}

\section{Concluding remarks}\label{S:6}
We have presented the computation of (i) multiplicative compensators; and (ii) additive and multiplicative compensators under a new measure obtained by the multiplicative compensation of a given non-negative semimartingale. In Theorems~\ref{T:mult comp} and \ref{T:Girsanov} we have treated these tasks as problems in their own right. 

In practice, the inputs to these computations are likely to come with more structure than indicated in the two theorems, i.e., the input processes $Z$ and $V$ in Theorems~\ref{T:mult comp} and \ref{T:Girsanov} will typically be represented with respect to some common underlying, possibly multivariate, process, say $X$. This is illustrated in Examples~\ref{E:190703.1}, \ref{E:201126}, \ref{E:MV}, and \ref{E:210410}, and Corollary~\ref{C:200616}, respectively. 

One should observe that here one may use the more general class of representing functions $\I(X)$ that are specific to $X$; see \cite[Definition~3.2]{cerny.ruf.22.ejp}. This is possible because a composition of a function in $\Uni$ with a compatible element of $\I(X)$ remains in $\I(X)$; see \cite[Corollary~3.20]{cerny.ruf.22.ejp}. For instance, this observation allows the use of general integrands $\zeta$ in Example~\ref{E:190703.1}.

\section*{Acknowledgements }
We thank Jan Kallsen, Johannes Muhle-Karbe, two anonymous referees, and an associate editor for helpful comments.\vspace{-0.0cm}

\def\MR#1{\href{http://www.ams.org/mathscinet-getitem?mr=#1}{MR#1}}
\def\ARXIV#1{\href{https://arxiv.org/abs/#1}{arXiv:#1}}
\def\DOI#1{\href{https://doi.org/#1}{doi:#1}}

\end{document}